\documentclass[11pt]{amsart}

\usepackage{hyperref}

\usepackage{enumerate}   
\usepackage{amsmath, amssymb, amsfonts, amsthm, fullpage}
\usepackage{graphics,graphicx}
\usepackage{stmaryrd}

\usepackage{marginnote}

\newcommand{\ind}{\mathrm{ind}}

\newcommand{\bo}{\mathbf{0}}

\newcommand{\bF}{\mathbf{F}}

\newcommand{\bu}{\mathbf{u}}

\newcommand{\bA}{\mathbf{A}}

\newcommand{\bG}{\mathbf{G}}

\newcommand{\bJ}{\mathbf{J}}
\newcommand{\bK}{\mathbf{K}}
\newcommand{\bU}{\mathbf{U}}

\newcommand{\bS}{\mathbf{S}}

\newcommand{\tbn}[1]{{\left\vert\kern-0.25ex\left\vert\kern-0.25ex\left\vert #1 \right\vert\kern-0.25ex\right\vert\kern-0.25ex\right\vert}}

\newcommand{\I}{\mathrm{i}}

\newtheorem{remark}{Remark}[section]

\newtheorem{proposition}{Proposition}[section]

\begin{document}

\title{Saddle-focus to saddle transition of 4D heteroclinic orbits resembling oscillatory and regularized shock waves of some dissipative Boussinesq systems}

\author{Rezwana Razzaque Angana}
\address{\textbf{R.~Angana:} Victoria University of Wellington, School of Mathematics and Statistics, PO Box 600, Wellington 6140, New Zealand}
\email{rezwana.angana@vuw.ac.nz}

\author{Brendan Harding}
\address{\textbf{B.~Harding:} Victoria University of Wellington, School of Mathematics and Statistics, PO Box 600, Wellington 6140, New Zealand}
\email{brendan.harding@vuw.ac.nz}

\author{Dimitrios Mitsotakis}
\address{\textbf{D.~Mitsotakis:} Victoria University of Wellington, School of Mathematics and Statistics, PO Box 600, Wellington 6140, New Zealand}
\email{dimitrios.mitsotakis@vuw.ac.nz}


\date{\today}


\begin{abstract}
Oscillatory and regularized shock waves have long been studied to understand the formation and dynamics of undular bores in shallow waters, phenomena that typically arise from the interplay between dispersion, nonlinearity, and dissipation. In this work, we investigate and describe the nature and shape of traveling wave solutions that resemble oscillatory and regularized shock waves of some Boussinesq-type water wave models known as Bona--Smith systems in the presence of dissipation. The associated traveling wave equations are solutions of four-dimensional dynamical systems, leading to both computational challenges and rich structural behavior. In particular, to describe the properties of these traveling waves, we analyze heteroclinic orbits that connect two distinct equilibria, one of which undergoes a transition in the four-dimensional phase space from a saddle-focus to a saddle point as the dissipation parameter increases. This transition governs the nature of the wave profile, which evolves smoothly from an oscillatory shock wave to a regularized (non-oscillatory) shock wave. To compute such solutions, we introduce a finite-difference-based numerical method tailored for this class of problems. 
\end{abstract}

\maketitle

\section{Introduction}

Undular bores are nonlinear wave phenomena that occur in shallow water when a rapid transition between two different flow states gives rise to an oscillatory wave train. Unlike classical hydraulic jumps, which result in a sharp discontinuity, undular bores are predominantly dispersive and exhibit a series of smooth, decaying oscillations. These structures typically emerge when an initial disturbance, such as a sudden change in water depth or velocity, propagates into a quiescent or slower-moving region. They are commonly observed in rivers and tidal flows and have been modeled effectively using weakly nonlinear dispersive equations \cite{Chanson2010, T1994}.

Nonlinear and dispersive water wave equations, such as the Korteweg--de Vries (KdV) and Benjamin--Bona--Mahony (BBM) equations, typically admit traveling wave solutions in the absence of dissipation, \cite{bbm,Whitham11}. These traveling waves correspond to homoclinic orbits (trajectories that leave and return to the same equilibrium) in the phase space of the associated dynamical systems, see for instance, \cite{Whitham11,W2003}. Such traveling water waves that can be represented by homoclinic orbits in the phase space are referred to as solitary waves, and they are often both stable and asymptotically stable, cf. \cite{B1972,PW1994}.

When the initial conditions of nonlinear dispersive equations (in the absence of dissipation) connect distinct asymptotic states, (i.e. different equilibria), the resulting dynamics often give rise to structures known as dispersive shock waves (DSWs) \cite{EH2016}. These feature an expanding region of rapid oscillations, reminiscent of undular bores, a resemblance that has prompted their use in modeling such physical phenomena \cite{Pere1966}. Unlike solitary waves, DSWs are not exact traveling wave solutions. Nevertheless, they display remarkable robustness and persist in the long-time dynamics of the system.

The formation of DSWs is often understood through their connection to traveling wave solutions in weakly dissipative perturbations of the original system. In such perturbed models, one typically observes the emergence of diffusive-dispersive shock waves, traveling wave solutions that arise due to the interplay between dispersion and weak dissipation \cite{EHS2017}. The asymptotic similarity between the original and perturbed models can explain why, in certain regimes, solutions originating from the same initial conditions tend to stay close over long times \cite{BMT2023,BMT2024}.

A well-studied example involves step-like initial conditions for the KdV equation that connect different asymptotic states. For small dissipation, the solutions of the KdV and KdV--Burgers equations remain close over finite time intervals, and the resulting DSWs can be viewed as perturbations of the diffusive-dispersive shock waves of the KdV--Burgers equation \cite{GP1987, EH2016, BS1985}. A similar behavior is expected for the BBM equation and its dissipative counterpart. However, the dynamics of the BBM Riemann problem are much richer than those of the KdV equation, and for larger jumps its solutions can differ substantially from those of the KdV and KdV--Burgers equations \cite{CEHS2021}. Similarly, for Boussinesq-type systems, such as the classical Boussinesq--Peregrine and Serre--Green--Naghdi systems, DSWs have also been shown to approximate the shock structures of related dissipative variants \cite{BMT2023,BMT2024}.

Diffusive-dispersive shock waves can be broadly classified into two types: oscillatory shock waves and regularized shock waves. Both connect distinct asymptotic states and are traveling wave solutions. An {\em oscillatory shock wave} is characterized by a sharp front followed by a wavetrain of decaying oscillations, and it typically arises in nonlinear dispersive equations with small but non-negligible dissipation. As the dissipation increases, the amplitude of the oscillations decreases, and beyond a certain threshold, the oscillatory structure disappears. The resulting traveling wave becomes strictly monotonic (decreasing if it is right-propagating, increasing if it is left-propagating), and is referred to as a {\em regularized shock wave}.

In this work, we investigate traveling wave solutions of certain  Boussinesq systems arising in water wave theory in the presence of dissipation, focusing in particular on a dissipative variant of the Bona--Smith systems \cite{BCS2002, BCS2004}, under the assumption of a flat bottom topography. The Bona--Smith systems are a family of coupled, nonlinear, dispersive partial differential equations that describe the bi-directional propagation of long-crested, small-amplitude surface water waves.

Let $\eta = \eta(x,t)$ denote the free surface elevation above the undisturbed water level at horizontal location $x \in \mathbb{R}$ and time $t \geq 0$, and let $u = u(x,t)$ be the horizontal velocity of the fluid measured at a fixed depth $z$. Denoting by $D > 0$ the undisturbed depth of the fluid, the dissipative Bona--Smith system takes the form:
\begin{equation}\label{eq:sysbs0}
\begin{aligned}
&\eta_t+Du_x+ (\eta u)_x- b D^2\eta_{xxt} = 0\ ,\\
&u_t+g\eta_x+ uu_x- cgD^2\eta_{xxx}- bD^2 u_{xxt} = \mu u_{xx}\ ,  
\end{aligned}
\end{equation}
where $g$ is the gravitational acceleration, and $\mu > 0$ is a diffusion parameter modeling the effect of fluid viscosity, excluding bottom friction \cite{DDZ2008, DD2007}. The parameters $b$ and $c$ are determined by the choice of the dimensionless constant $\theta^2$ (with $\frac{2}{3} \leq \theta^2 \leq 1$) through the relations
$$b=\frac{3\theta^2-1}{6} \qquad \text{and}\qquad c=\frac{3\theta^2-2}{3}\ .$$
This system, supplemented with appropriate initial conditions $\eta(x,0) = H_0(x)$ and $u(x,0) = U_0(x)$, has been shown to be globally well-posed in suitable Sobolev spaces \cite{Hu2009}.

The dissipative term $\mu u_{xx}$ in (\ref{eq:sysbs0}) models the dissipation due to the viscosity of the fluid. Boussinesq systems with dissipative terms of this form were derived from the Navier--Stokes equations using the theory of weakly damped potential flows in \cite{DDZ2008, DD2007}, where the viscous dissipation appears only in the momentum equation. A similar derivation of the dissipative Boussinesq--Peregrine system from the Euler equations with dissipation, in which the horizontal component of the viscous diffusion is neglected, was presented in \cite{BMT2022}. Since the Bona--Smith systems are asymptotically equivalent to the Boussinesq--Peregrine system, the same arguments justify the dissipative term in (\ref{eq:sysbs0}) and similar derivation can be constructed. Moreover, $u_{xx}$ is the simplest dissipative term that turns the Bona--Smith systems into diffusive-dispersive regularizations of the shallow water equations, in analogy with the KdV--Burgers equation.

Unlike the KdV and BBM equations, which describe unidirectional wave propagation, the Bona--Smith systems describe wave motion in both spatial directions. In the absence of dissipation, these systems are known to admit solitary wave solutions as traveling waves \cite{Chen1998, Chen2000, DM2008}. However, to the best of our knowledge, no results have been reported related to traveling waves in the presence of dissipation. For this reason, we study  traveling wave solutions of the dissipative Bona--Smith systems that correspond to heteroclinic orbits connecting two distinct equilibria of the corresponding dynamical systems. 

Beginning with a theoretical analysis of the stability properties of equilibria in the associated four-dimensional dynamical systems, we establish the existence of a critical value $\mu_0$ of the dissipation parameter $\mu$ at which a bifurcation occurs, corresponding to a transition from a saddle-focus to a saddle equilibrium. This result generalizes the well-known focus-to-node bifurcation observed in lower-dimensional settings, particularly in two-dimensional dynamical systems arising in other Boussinesq-type models, such as the classical Boussinesq--Peregrine system, and the KdV--Burgers equation \cite{BMT2022, BS1985}.

It is worth noting that studies addressing saddle-focus transitions, especially in four-dimensional systems, are relatively rare. Existing literature typically concentrates on homoclinic bifurcations associated with such equilibria \cite{KKM2019, Kuz2024}.

The existence of a critical value $\mu_0$ for the dissipation parameter, and the associated transition of the equilibrium from a saddle-focus to a saddle, have a direct impact on the structure of the heteroclinic orbits. As in the two-dimensional case, the traveling wave solutions transition from oscillatory shock waves for $\mu < \mu_0$ to regularized (strictly monotonic) shock waves for $\mu \geq \mu_0$. However, since the underlying dynamical systems are four-dimensional, the existence of additional heteroclinic orbits with qualitatively different profiles cannot be ruled out.

These findings are based on numerical computations, which indicate that the heteroclinic orbits are structurally stable and form a continuous family of solutions with respect to the parameter $\mu$. Furthermore, these orbits appear to be natural continuations of symmetric homoclinic orbits corresponding to the undamped case $\mu = 0$. This scenario is characteristic of homoclinic bifurcation phenomena and reflects the transverse intersection of the stable and unstable manifolds of the equilibria connected by the heteroclinic orbit.

Capturing heteroclinic orbits in four-dimensional dynamical systems presents notable challenges, primarily due to their non-uniqueness. These orbits are inherently non-unique not only because of their invariance under horizontal translations, but also because additional, potentially nonphysical, heteroclinic connections may arise from homoclinic bifurcations involving multi-pulse homoclinic orbits, as described in \cite{Chen2000}.

To address these difficulties, we introduce a finite difference method specifically tailored to compute heteroclinic orbits that resemble oscillatory and regularized shock waves. The method is applied to the corresponding four-dimensional ordinary differential equations, with boundary conditions imposed at the two equilibria connected by the heteroclinic orbit. Appropriate asymptotic conditions are enforced to promote convergence and suppress nonphysical solutions. To circumvent the non-uniqueness arising from translational invariance, the solution is centered so that the velocity component $u$ has an inflection point at $x = 0$, providing a consistent and physically meaningful reference frame.

The structure of the paper is as follows. Section \ref{sec:review} provides an overview of the basic properties of the dissipative Bona--Smith systems, along with a brief discussion of the Peregrine (classical Boussinesq) system. In Section \ref{sec:existence}, we focus (primarily for simplicity) on the BBM--BBM system with $c = 0$, and study the corresponding dynamical system arising from the traveling wave ansatz. In particular, we analyze the dependence of the equilibria on the parameter $\mu$, revealing a saddle-focus to saddle bifurcation at a critical value $\mu_0$. This bifurcation plays a key role in determining the qualitative structure of the traveling waves. We also discuss some properties of the associated heteroclinic orbits. In Section \ref{sec:bonasmith}, we extend the analysis to the full Bona--Smith system with $c > 0$. Although the resulting expressions are more involved, the qualitative behavior of the underlying dynamical system remains consistent with the case $c = 0$.
Section \ref{sec:computho} introduces a numerical method designed to compute heteroclinic orbits that resemble undular bores.  Section \ref{sec:comparison} presents a comparison of traveling wave profiles across different systems, computed using the numerical method introduced in Section \ref{sec:computho}. This includes an exploration of how the parameter $\mu$ influences the shape of the resulting waves. Section \ref{sec:conclusions} summarizes the main findings of the study. Finally, the detailed analysis of the eigenvalues of the linearization about the nontrivial equilibrium, which leads to the saddle-focus to saddle transition stated in Section \ref{sec:existence}, is presented in Appendix \ref{app:linearization}.

\section{Review of the mathematical models}\label{sec:review}

The dissipative Bona--Smith family of systems (\ref{eq:sysbs0}) in dimensionless and scaled form can be written as
\begin{equation}\label{eq:sysbs1}
\begin{aligned}
&\eta_t+u_x+ \varepsilon (\eta u)_x-\sigma^2 b \eta_{xxt} = 0\ ,\\
&u_t+\eta_x+ \varepsilon uu_x- \sigma^2c\eta_{xxx}-\sigma^2 b u_{xxt} = \nu u_{xx}\ ,  
\end{aligned}
\end{equation}
where for $\frac{2}{3}\leq \theta^2\leq 1$ the coefficients $b,c$ are related via the formulas
$$b=\frac{3\theta^2-1}{6} \qquad \text{and}\qquad c=\frac{3\theta^2-2}{3}\ .$$ 
The parameters $0 < \varepsilon < 1$ and $0 < \sigma < 1$ characterize, respectively, the ratio of wave amplitude to undisturbed water depth and the ratio of water depth to wavelength \cite{BCS2002}. In the scaled variables the horizontal bottom is at $D=-1$ and the scaled diffusion parameter is given by $\nu = \varepsilon \kappa > 0$, where $\kappa = O(1)$ describes dissipative effects due to internal viscosity.

The case $\theta^2 = 2/3$ corresponds to the so-called BBM--BBM system, which will be considered separately due to its relative simplicity. However, as we shall demonstrate later, its behavior closely resembles that of the more general Bona--Smith systems with $2/3 < \theta^2 \leq 1$. The system with $\theta^2 = 1$ is commonly referred to as the classical Bona--Smith system, as it was the original model derived by Bona and Smith in their foundational work \cite{BS}.

Boussinesq systems provide valid approximations to the full Euler equations in the long-wave regime when the Stokes number $S = \varepsilon / \sigma^2$ is of order one, i.e., $S = O(1)$ \cite{Whitham11}. For simplicity, we assume $\varepsilon = \sigma^2$, which satisfies this condition.

Applying the change of variables:
\begin{equation}\label{eq:chvar}
\eta \to \varepsilon \eta, \quad u \to \varepsilon u, \quad x \to \frac{x}{\sqrt{\varepsilon b}}, \quad t \to \frac{t}{\sqrt{\varepsilon b}}\ ,
\end{equation}
we rescale system \eqref{eq:sysbs1} to obtain the dimensionless form:
\begin{equation}\label{eq:sysbs2}
\begin{aligned}
&\eta_t+u_x+ (\eta u)_x- \eta_{xxt} = 0\ ,\\
&u_t+\eta_x+ uu_x- \gamma \eta_{xxx}- u_{xxt} = \delta u_{xx}\ , 
\end{aligned}
\end{equation}
where
$$\gamma=\frac{c}{b}=\frac{2(3\theta^2-2)}{3\theta^2-1} \qquad\text{and}\qquad \delta=\frac{\nu}{\sqrt{b\varepsilon}}  \ .$$
Note that for $\theta^2\in [2/3,1]$, the parameter $\gamma$ lies in the interval $[0,1]$.

Consider the impulse functional
\begin{equation}\label{eq:impulse}
I(t) = \int_{-\infty}^{+\infty} \left( u \eta + u_x \eta_x \right) ~ dx\ .
\end{equation}
Then, for solutions of system \eqref{eq:sysbs2} for which $I(t)$ is well-defined, the time derivative satisfies
\begin{equation}\label{eq:pdedis}
\frac{d}{dt} I(t) = -\delta \int_{-\infty}^{+\infty} u_x \eta_x ~ dx\ .
\end{equation}
Similarly, if
\begin{equation}
E(t) = \int_{-\infty}^{+\infty} \left( \eta^2 + \gamma \eta_x^2 + (1 + \eta) u^2 \right) ~ dx
\end{equation}
is the energy functional, then
\begin{equation}
\frac{d}{dt} E(t) = 2\delta \int_{-\infty}^{+\infty} u_x \eta_t ~  dx\ .
\end{equation}
These relations do not constitute dissipation laws in the usual sense. On the other hand, solutions of system \eqref{eq:sysbs2} for which the previous integrals make sense have been studied in \cite{CG2007}, where it was shown that their $L^2$-norms decay at the rate $O(t^{-1/4})$. Consequently, the Bona--Smith-type systems (\ref{eq:sysbs2}) with $2/3 \leq \theta^2 \leq 1$ can be viewed as dissipative systems. Moreover, for systems with $\theta^2 > 2/3$, it was also proved in \cite{CG2007} that the $H^1$-norm of $\eta$ decays in time at the same rate, $O(t^{-1/4})$. However, the existing theoretical framework does not cover solutions with unbounded norms, and particularly those corresponding to undular bores that do not decay to zero as $x \to -\infty$. In this work, we focus on such positive traveling wave solutions, characterized by the asymptotic behavior
$$\lim_{x\to-\infty}(u,\eta)=(u_0,\eta_0)\qquad \text{and} \qquad \lim_{x\to+\infty}(u,\eta)=(0,0)\ ,$$
where $u_0>0$ and $\eta_0>0$.

An asymptotically equivalent system that has been employed for the study of undular bores is the classical Boussinesq--Peregrine system \cite{Pere1966,Pere1967,BMT2022} and in scaled form is written as
\begin{equation}\label{eq:cbscaled}
    \begin{aligned}
&\eta_t+u_x+ \varepsilon(\eta u)_x = 0\ ,\\
&u_t+\eta_x+ \varepsilon uu_x- \sigma^2 b u_{xxt} = \nu u_{xx}\ .
    \end{aligned}
\end{equation}
Such systems were analyzed in \cite{AAG2007}.

Using the same change of variables as in the case of the Bona--Smith systems, we obtain the system
\begin{equation}\label{eq:cbscaled2}
\begin{aligned}
&\eta_t + u_x + (\eta u)_x = 0\ ,\\
&u_t + \eta_x + uu_x - u_{xxt} = \bar{\delta} u_{xx}\ ,
\end{aligned}
\end{equation}
where $\bar{\delta} = \nu/\sqrt{b\varepsilon}$ with $b=1/3$ for the classical Boussinesq--Peregrine system. For the same values of $\nu$ and $\varepsilon$, this parameter is related to the parameter $\delta$ of the Bona--Smith systems (\ref{eq:sysbs2}) via $\bar{\delta}=\delta/\sqrt{2-\gamma}$. In particular, $\bar{\delta}=\delta$ for $\gamma=1$ and $\bar{\delta}=\delta/\sqrt{2}$ for the BBM--BBM system. Traveling waves of this particular system were studied in \cite{BMT2022}, while its dissipative nature was analyzed in \cite{AAG2007,CG2007}.

System (\ref{eq:cbscaled2}) admits classical solitary wave solutions in the form of traveling waves for all wave speeds $s > 1$ when $\bar{\delta} = 0$ \cite{Chen2000,DM2008}. These solutions correspond to symmetric homoclinic orbits to the origin in a two-dimensional phase space and undergo a homoclinic--heteroclinic bifurcation as $\bar{\delta}$ transitions from zero to positive values \cite{BMT2022}. More precisely, for all $s > 1$ and $\bar{\delta} > 0$, there exists a unique heteroclinic orbit connecting the origin to a positive equilibrium, given by
$$    u_0=\frac{3s-\sqrt{s^2+8}}{2}\quad\text{and}\quad \eta_0=\frac{u_0}{s-u_0}\ .$$

It was further shown in \cite{BMT2022} that there exists a critical bifurcation value $\bar{\delta}_0$, which can be computed analytically as
\begin{equation}
\bar{\delta}_0 = \sqrt{4s\,\alpha(s)}, \qquad \text{where} \quad \alpha(s) = \frac{s - \sqrt{s^2 + 8}}{2} + \frac{4s}{(s - \sqrt{s^2 + 8})^2}\ ,
\end{equation}
such that the nature of the heteroclinic orbits depends on the value of $\bar{\delta}$. For $\bar{\delta} < \bar{\delta}_0$, the solutions take the form of oscillatory shock waves, while for $\bar{\delta} \geq \bar{\delta}_0$, the resulting waves are regularized shock waves, i.e. smooth, monotonic traveling waves without oscillations. When the dispersion parameter $b$ is diminishing at a particular rate with respect to an also diminishing diffusion parameter $\nu$, these traveling waves converge to entropic shock waves of the non-dispersive shallow water equations \cite{BMT2023}.

Although systems \eqref{eq:sysbs1} and \eqref{eq:cbscaled} are asymptotically equivalent in the long-wave, weakly nonlinear regime, they differ fundamentally in their mathematical structure. In particular, the phase space associated with traveling wave solutions of the Bona--Smith systems is four-dimensional, while for the classical Boussinesq--Peregrine system it is only two-dimensional. This key distinction raises natural questions about the influence of dimensionality on the underlying dynamics and how the added complexity affects the qualitative features of the corresponding wave solutions.

\section{Dynamics of heteroclinic orbits of dissipative BBM--BBM system}\label{sec:existence}

For the sake of clarity, we examine in detail the traveling wave dynamics of the dissipative BBM--BBM system ($\gamma=0$). The more general dissipative Bona--Smith system with $\gamma \in (0,1]$ is considered only briefly, since its analysis is more involved and several properties of the BBM--BBM system do not extend directly to this case.

For $\delta>0$, consider the BBM--BBM system 
\begin{equation}\label{eq:bbmbbm}
    \begin{aligned}
\eta_t+u_x+ (\eta u)_x- \eta_{xxt} &= 0\ ,\\
u_t+\eta_x+ uu_x- u_{xxt} &= \delta u_{xx}\ .
    \end{aligned}
\end{equation}
We search for right-traveling waves for $s>1$ in the form:
\begin{equation}\label{eq:ansatz}
    u(\xi)=u(x-s t)>0,\quad \eta(\xi)=\eta(x- s t)>0\ ,
\end{equation}
with the boundary conditions
\begin{equation}\label{eq:bc1}
    \lim_{\xi\to +\infty}(\partial_\xi^{j} u(\xi),\partial_\xi^j \eta(\xi))=(0,0)\quad \text{for $j=0,1,2,\dots$}\ ,
\end{equation}
while 
\begin{equation}\label{eq:bc2}
    \begin{aligned}
      &\lim_{\xi\to -\infty}(u(\xi),\eta(\xi))=(u_0,\eta_0)\ , \\
      &\lim_{\xi\to -\infty }(\partial_\xi^{j} u(\xi),\partial_\xi^j \eta(\xi))=(0,0) \quad \text{for $j=1,2,\dots$}\ .
    \end{aligned}
\end{equation}

Such special solutions represent heteroclinic orbits connecting two distinct equilibria of the system 
\begin{equation}\label{eq:sys1}
\begin{aligned}
 -s \eta''+s\eta-u-\eta u & =0\ ,\\
 -s u''+su-\eta-\tfrac{1}{2} u^2 &= -\delta u'\ ,
\end{aligned}
\end{equation}
which is obtained after substituting the {\em ansatz} (\ref{eq:ansatz}) into (\ref{eq:bbmbbm}), and integrating over $(\xi,+\infty)$ for some $\xi\in\mathbb{R}$.

Physically relevant traveling-wave solutions of system (\ref{eq:bbmbbm}) are expected to exist for wave speeds satisfying $|s| > 1$, since $|s| = 1$ corresponds to the limiting case of the linear, infinitely long wave speed of propagation.
The system (\ref{eq:sys1}) possesses the symmetry
\begin{equation}\label{eq:rev1}
\xi \to -\xi,\quad u \to -u,\quad \eta \to \eta,\quad s \to -s\ ,
\end{equation}
which implies that if $(u(\xi), \eta(\xi), s)$ is a traveling-wave solution of (\ref{eq:bbmbbm}) with speed $s > 1$, then the reflected profile $(-u(\zeta), \eta(\zeta), -s)$, with $\zeta = -\xi$, is also a solution, now with speed $s < -1$.
Owing to this symmetry, it suffices to consider traveling waves with $s > 1$. Such waves propagate to the right, whereas those with $s < -1$ propagate to the left.

\subsection{Dynamical picture}\label{sec:dynamics}

System (\ref{eq:sys1}) can be written in the form
\begin{equation}\label{eq:sys2}
    \bS \bu'' + \nabla f(\bu) =  \delta~ \bK \bu'\ ,
\end{equation}
where $\bu=(u,\eta)^T$, $\bu'=(u',\eta')^T$,
\begin{equation}
\bS=\begin{pmatrix}
    0 & -s\\
    -s & 0
\end{pmatrix}, \quad \bK=\begin{pmatrix}
    0 & 0\\
    -1 & 0
\end{pmatrix} \quad \text{and} \quad
f(\bu)=s~u\eta  -\tfrac{1}{2}(\eta^2+u^2+ u^2 \eta)\ .
\end{equation}

\begin{figure}[t]
\centering
\includegraphics[width=0.8\columnwidth]{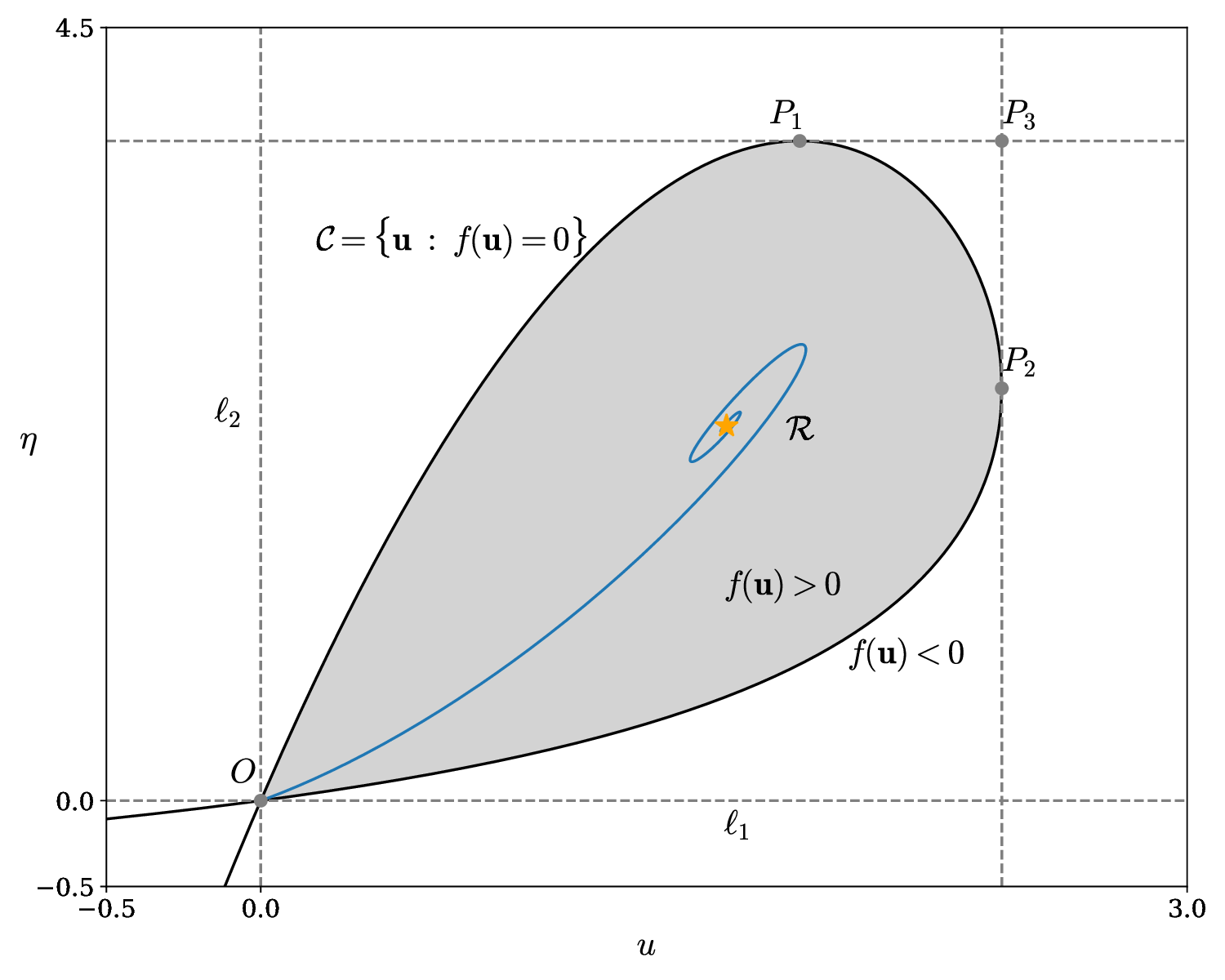}
\caption{A projection of a heteroclinic orbit connecting the critical points $\bo$ and $\bU_0$ (depicted by a star) on physical space when $s=2.2$ and $\delta=1.8$. The projected orbit $\mathcal{R}$ lies entirely inside the interior of the shaded domain $\mathcal{D}$.}
\label{fig:domain}
\end{figure}

The special case $\delta = 0$ has been previously studied in \cite{Toland1986, Chen2000}, where the system (\ref{eq:sys1}) preserves the Hamiltonian
$$H(\bu)=Q(\bu')+f(\bu)\ ,$$
with $Q(\bu) = \tfrac{1}{2} \bu^T \bS \bu = -su\eta$ representing an indefinite quadratic form (i.e., $\det(\bS) < 0$), with $Q(\bu)<0$ in the first and third quadrants and $Q(\bu)>0$ in the second and fourth quadrants of the plane. Moreover, we have $Q(\bu)=0$ along the lines $\ell_1$ and $\ell_2$, which in the BBM--BBM case coincide with the $u$-axis and the $\eta$-axis, respectively. The function $f$ satisfies $f(0, 0) = 0$, and its zero-level set 
$$\mathcal{C}=\{\bu\in\mathbb{R}^2:~ f(\bu)=0\}\ ,$$ 
lies in the region where $Q(\bu) \leq 0$, namely the closed first and third quadrants of the physical $(u, \eta)$-plane. Its part in the first quadrant is a closed convex curve, which is smooth except at the origin, where it has a corner with tangent lines $\eta=(s\pm\sqrt{s^2-1})u$. Its part in the third quadrant consists of two unbounded branches. Typically, this curve has an inverted $\gamma$-like shape in the $(u, \eta)$-plane. Figure \ref{fig:domain} illustrates portions of $\mathcal{C}$ in the first and third quadrants. Notably, the tangent to $\mathcal{C}$ at the point $P_1 = ((s^2 - 1)/s, s^2 - 1)$ is horizontal (i.e., parallel to $\ell_1$), while at $P_2 = (2(s - 1), 2(s - 1))$ it is vertical (i.e., parallel to $\ell_2$), i.e. they are parallel to the lines at which $Q(\bu)=0$. The parallel lines to $\ell_1$ and $\ell_2$ passing through $P_1$ and $P_2$ respectively, intersect at the point $P_3=(2(s-1),s^2-1)$.

In the following, we refer to by $\mathcal{R} = \{\bu(\xi) = (u(\xi), \eta(\xi)) : \xi \in \mathbb{R}\}$ any orbit of system (\ref{eq:sys2}), where $\bu = (u, \eta)$. The domain enclosed by $\mathcal{C}$ in the first quadrant is bounded by the arcs $\Gamma_1$ and $\Gamma_2$, which connect the origin $O$ to the points $P_1$ and $P_2$, respectively, along with the arc $P_1P_2$ forming a convex and closed region.

Although we do not provide a theoretical justification, numerical evidence presented later suggests that heteroclinic orbits of the Bona--Smith systems (including the BBM--BBM system) that exhibit the structure of oscillatory or regularized shock waves lie entirely within the interior of the domain
$$\mathcal{D}=\mathring{\mathcal{C}}\cap\{u>0,~\eta>0\}\ ,$$ 
where $\mathring{\mathcal{C}}$ denotes the interior of the curve $\mathcal{C}$. If an intersection between $\mathcal{R}$ and $\mathcal{C}$ can occur, this seem to be on a very small scale along the arc $P_1P_2$ and only for very small values of the parameter $\delta$. This region is represented as the shaded area in Figure \ref{fig:domain}. Specifically, one can show that for heteroclinic orbits of (\ref{eq:sys2}) with $H(\bU)>0$ for all $\xi\in\mathbb{R}$, are confined in the region $\mathring{\mathcal{C}}\cup\mathcal{T}$, where $\mathcal{T}$ is the area bounded by the arcs $P_1P_2$, $P_1P_3$ and $P_2P_3$, including the open arc $P_1P_2$. It is noted that in the absence of dissipation, symmetric homoclinic orbits to the origin lie inside $\mathcal{D}$ but they achieve their maximum value on the arc $P_1P_2$ \cite{Toland1986}.

Solving the equation $f(\bu)=0$ with respect to $u$ we obtain \cite{Chen2000}:
$$u_\pm=\frac{s\eta\pm \eta\sqrt{s^2-1-\eta}}{1+\eta}\ ,$$ 
which is meaningful when $0<\eta\leq s^2-1$. Similarly, solving with respect to $\eta$ we obtain:
$$\eta_\pm=u\frac{2s - u \pm \sqrt{(u-2(s+1))(u-2(s-1))}}{2}\ ,$$
which is meaningful when $0<u\leq 2(s-1)$. Clearly, 
$$f(\bu)=-\tfrac{1}{2}(1+\eta)(u-u_-)(u-u_+)\ ,$$
and $f(\bu)>0$ for all $\bu\in \mathring{\mathcal{C}}$. Moreover, $\nabla f(\bu)\not=\bo$ for all $\bu\in\mathcal{C}$ except at $\bu=\bo$. An example of the curve $\mathcal{C}$ and an embraced heteroclinic orbit $\mathcal{R}$ is depicted in Figure \ref{fig:domain}.

The dynamical system (\ref{eq:sys2}) admits three equilibrium points, corresponding to the critical points of the function $f$, i.e., where $\nabla f(\bu) = 0$. These equilibria are determined by solving the system
\begin{equation}\label{eq:critic}
su-\frac{u}{s-u}-\frac{1}{2}u^2=\frac{u(u^2-3s u+2(s^2-1))}{2(s-u)}=0\quad \text{and}\quad \eta=\frac{u}{s-u}\ .
\end{equation}
The solutions of this system are the origin $\bo = (0, 0)$ and the two nontrivial equilibria $\bu_0 = (u_0, \eta_0)$ and $\bu_1 = (u_1, \eta_1)$, where
\begin{equation}\label{eq:equilibria}
    u_0=\frac{3s-\sqrt{s^2+8}}{2}\quad\text{and}\quad \eta_0=\frac{u_0}{s-u_0}=-1 +\tfrac{1}{4}(s^2 + s \sqrt{s^2+8})\ ,
\end{equation}
and 
$$u_1=\frac{3s+\sqrt{s^2+8}}{2}\quad \text{and}\quad \eta_1=\frac{u_1}{s-u_1}=-1 +\tfrac{1}{4}(s^2 - s \sqrt{s^2+8}) \ .$$
For $s > 1$, the equilibrium $\bu_1$ lies in a nonphysical region, where the wave goes below the impermeable bottom since $\eta_1 < -1$, and is therefore excluded from the present analysis. Note that the quadratic polynomial $u^2-3su+2(s^2-1)$ in (\ref{eq:critic}) is positive for $u<u_0$ and $u>u_1$.

We note that the critical point $(u_0,\eta_0)\in\mathring{\mathcal{C}}$ since \begin{equation}\label{eq:fuz}
f(\bu_0)=\left(\tfrac{3s-\sqrt{s^2+8}}{2(s-\sqrt{s^2+8})}\right)^2(-2-s^2+s\sqrt{s^2+8})>0 \quad \text{for all $s>1$}\ .
\end{equation}
Moreover, $f(\bu_0)$ is a local maximum of the function $f(\bu)$ within $\mathring{\mathcal{C}}$. Indeed, $\nabla f(\bu_0)=0$, and the second order conditions yield $f_{uu}(\bu_0)=-(1+\eta_0)<0$ and 
$$ D^2 f(\bu_0)=
\begin{vmatrix}
f_{uu}(\bu_0) & f_{u\eta}(\bu_0)\\
f_{\eta u}(\bu_0) & f_{\eta\eta}(\bu_0)
\end{vmatrix}=1+\eta_0-(s-u_0)^2>0\ ,
$$
for all $s>1$. (The last inequality holds because $u_0>s-\sqrt[3]{s}$). In contrast, at the origin $\bo=(0,0)$ we have $f_{uu}(\bo)=-1<0$ and $D^2f(\bo)=1-s^2<0$, for $s>1$, implying that the origin is a saddle point (a local minimum within the closure $\bar{\mathcal{D}}$). The point $\bu_1$ is also a saddle point of $f(\bu)$, since $D^2f(\bu_1)=1+\eta_1-(s-u_1)^2<0$ because $\eta_1<-1$. Since there are no other critical points in $\mathcal{D}$, we conclude that any solution of (\ref{eq:sys2}) that lies in the interior $\mathcal{D}$ satisfies $0<f(\bu)\leq f(\bu_0)$. 

To analyze system (\ref{eq:sys1}) as a four-dimensional dynamical system, we recast it in the form
\begin{equation}\label{eq:sys4d}
\bU'=\bF(\bU)\doteq\bA\bU+\bG(\bU)\ ,
\end{equation}
where $\bU=(u,\eta,v,h)^T$, 
\begin{equation}
\bA=\begin{pmatrix}
    0 & 0 & \tfrac{1}{s} & 0\\
    0 & 0 & 0 & \tfrac{1}{s} \\
    s & -1 & \tfrac{\delta}{s} & 0\\
    -1 & s & 0 & 0
\end{pmatrix} \quad \text{and}\quad \bG(\bU) = \begin{pmatrix}
    0\\
    0\\
    -\frac{1}{2}u^2\\
    -u\eta 
\end{pmatrix}\ .
\end{equation}
The physical space $(u, \eta)$ of system (\ref{eq:sys2}) corresponds to the natural projection of the full four-dimensional phase space $(u, \eta, v, h)$ associated with system (\ref{eq:sys4d}) onto the $(u, \eta)$-plane.
System (\ref{eq:sys4d}) possesses two physically relevant equilibrium points: the origin $\bo = (0, 0, 0, 0)^T$, and $\bU_0 = (u_0, \eta_0, 0, 0)^T$, where $u_0$ and $\eta_0$ are given in (\ref{eq:equilibria}). Linearizing around the origin yields a linear system governed by the matrix $\bA$, whose eigenvalues are the roots of the characteristic polynomial 
\begin{equation}\label{eq:P1poly}
P_1(\lambda)=\lambda^4 -\frac{\delta}{s}\lambda^3- 2 \lambda^2 + \frac{\delta}{s}\lambda + \left(1-\frac{1}{s^2}\right)\ .    
\end{equation}
Applying Descartes' Rule of Signs (see \cite{W2003}) to $P_1(\lambda)$ shows that there are at most two positive real roots. Similarly, applying the same  rule to $P_1(-\lambda)$ shows that there are at most two negative real roots. To determine the exact count, we observe that $P_1(0)=1-1/s^2>0$, $P_1(\pm 1)=-1/s^2<0$ and $P_1(\pm\infty)=+\infty$. These sign changes imply that the polynomial has exactly two positive and two negative real roots. Therefore, for all $\delta \geq 0$ and $s > 1$, the origin is a saddle point with two unstable and two stable directions in phase space. Moreover, the eigenvalues can be classified by magnitude: the leading eigenvalues lie in the interval $(-1, 1)$, while the remaining two eigenvalues lie in $(1, +\infty)$ and $(-\infty, -1)$, respectively.

When $\delta = 0$, the matrix $\bA$ admits the explicit eigenvalues
$$\rho_1=-\sqrt{1+\frac{1}{s}}<-1,\quad \rho^{-}=-\sqrt{1-\frac{1}{s}}>-1,\quad \rho^{+}=\sqrt{1-\frac{1}{s}}<1,\quad 
\rho_2=\sqrt{1+\frac{1}{s}}>1\ ,$$
where we denote by $\rho^{-}$ and $\rho^{+}$ the leading eigenvalues, i.e., those closest to the origin. In this case, the saddle quantity, which is defined as the sum of the leading stable and unstable eigenvalues, is given by $\sigma_0=\rho^{+}+\rho^{-}=0$. This is a degenerate situation, as $\sigma_0 = 0$ lies precisely on the boundary between dissipative and conservative saddle behavior. To the best of our knowledge, such a case is not covered by the standard theory of homoclinic bifurcations; see, for instance, \cite{Kuz2024}. Throughout this work, we use the notation $\rho$ to denote the eigenvalues of an equilibrium point when $\delta = 0$, and $\lambda$ when $\delta > 0$.

Furthermore, assuming that for small $\delta > 0$, the eigenvalues of $\bA$ admit an asymptotic expansion of the form $\lambda=\lambda^0+\delta \lambda^1+O(\delta^2)$, we observe that the roots of the characteristic polynomial (\ref{eq:P1poly}) shift to the right as $\delta$ increases from zero. In particular, for sufficiently small $\delta$, each eigenvalue is perturbed from its $\delta = 0$ value by the same leading-order term, namely $\delta / (4s)$. That is, as $\delta \to 0$, we have the expansions:
$$\lambda_1\sim \frac{\delta}{4 s}+\rho_1,\quad \lambda_{-}\sim\frac{\delta}{4 s}+\rho_{-},\quad \lambda_{+}\sim\frac{\delta}{4 s}+\rho_{+},\quad \lambda_2\sim \frac{\delta}{4 s}+\rho_2\ .$$
It follows that the saddle quantity for $\delta > 0$, defined as $\sigma_\delta = \lambda_{+} + \lambda_{-}$, satisfies the asymptotic relation $\sigma_\delta\sim \delta/2s >0$, indicating a shift from the degenerate case $\sigma_0 = 0$ (for $\delta = 0$) to a non-degenerate saddle for $\delta > 0$. Finally, as $\delta \to \infty$, the eigenvalues retain their qualitative ordering within their respective intervals, and their asymptotic behavior is given by
$$\lambda_1\sim -1 - \frac{1}{2 \delta s},\quad \lambda_{-}\sim\frac{1-s^2}{\delta s},\quad \lambda_+\sim 1-\frac{1}{2 \delta s},\quad \lambda_2\sim \frac{\delta}{s}\ .$$
Notably, despite the dissipation limit, the saddle quantity remains positive, preserving the system's saddle-type nature.

To analyze the nature of the equilibrium $\bU_0$, we examine the linearization matrix
\begin{equation}\label{eq:linatild}
\tilde{\bA}=\tilde{\bA}(\delta)=\begin{pmatrix}
0& 0& \tfrac{1}{s}& 0\\
0& 0& 0& \tfrac{1}{s}\\
s - u_0& -1& \tfrac{\delta}{s} &0 \\ 
-(1 + \eta_0)& s - u_0& 0& 0
\end{pmatrix}\ .
\end{equation}
Its eigenvalues satisfy the characteristic equation
\begin{equation}\label{eq:P2poly}
P_2(\lambda)=\lambda^4-\frac{\delta}{s}\lambda^3-2\frac{s-u_0}{s}\lambda^2+\frac{\delta}{s}\frac{s-u_0}{s}\lambda-\frac{1+\eta_0-(s-u_0)^2}{s^2}\ .
\end{equation}
For the coefficients of $P_2$ involving $u_0$ and $\eta_0$ and for $s>1$, we have that
$$A\doteq\frac{s-u_0}{s}=\frac{-s+\sqrt{s^2+8}}{2s}>0\ ,$$
and 
\begin{equation}\label{eq:q}
q\doteq\frac{1+\eta_0-(s-u_0)^2}{s^2}=\frac{-(s^2+8)+3s\sqrt{s^2+8}}{4s^2}=\frac{(3s-\sqrt{s^2+8})\sqrt{s^2+8}}{4s^2}>0\ .
\end{equation}

When $\delta=0$, the polynomial $P_2(\lambda)$ has the four eigenvalues $\rho_i^\pm=\pm \I \rho_i$ and $\rho_r^\pm=\pm \rho_r$, where
\begin{equation}\label{eq:rootsdelta0}
        \rho_i=\sqrt{\frac{u_0-s+\sqrt{1+\eta_0}}{s}},\quad \rho_r=\sqrt{\frac{s-u_0+\sqrt{1+\eta_0}}{s}}\ .
\end{equation}
It is noted that $\frac{u_0-s+\sqrt{1+\eta_0}}{s}>0$ and $\frac{s-u_0+\sqrt{1+\eta_0}}{s}>0$ for all $s>1$.

In what follows, we will denote the real roots of $P_2(\lambda)$ by $\lambda_r^+>0$, $\lambda_r^-<0$  and $\lambda_i^+$ and $\lambda_i^-$ the roots of $P_2(\lambda)$ for $\delta>0$ that emerge from the imaginary roots of $P_2(\lambda)$ for $\delta=0$ as $\delta$ is increasing.

The following proposition summarizes the behavior of the eigenvalues of $\tilde{\bA}(\delta)$ as $\delta$ varies. Its proof relies on standard but lengthy computations and is given in Appendix \ref{app:linearization}.

\begin{proposition}\label{prop:eigenvalues}
Let $s>1$. There is a unique critical value $\delta_0=\delta_0(s)>0$, given explicitly by (\ref{eq:critdelta0}), such that the following hold:
\begin{enumerate}[(i)]
\item For all $\delta>0$, the roots $\lambda_r^\pm$ are real and simple, $\lambda_r^-$ is the only negative root of $P_2(\lambda)$ and $\lambda_r^+$ is its largest real root. Both are increasing functions of $\delta$ and satisfy $\lambda_r^-<-\sqrt{A}$, $\lambda_r^+>\sqrt{A}$ and $\lambda_r^++\lambda_r^->0$.
\item For $0<\delta<\delta_0$, the roots $\lambda_i^\pm$ are non-real and complex conjugate with $0<Re(\lambda_i^\pm)<\delta/(2s)$ and $Re(\lambda_i^\pm)<\lambda_r^+$. Moreover, $Im(\lambda_i^+)$ is strictly decreasing in $\delta$.
\item For $\delta=\delta_0$, the roots $\lambda_i^\pm$ coalesce into a double positive root $\lambda_d\in(0,\sqrt{A})$.
\item For $\delta>\delta_0$, the roots $\lambda_i^\pm$ are real with $0<\lambda_i^-<\lambda_i^+<\sqrt{A}$. Moreover, $\lambda_i^-$ is decreasing and $\lambda_i^+$ is increasing in $\delta$.
\end{enumerate}
\end{proposition}

\begin{figure}[t]
\centering
\includegraphics[width=\columnwidth]{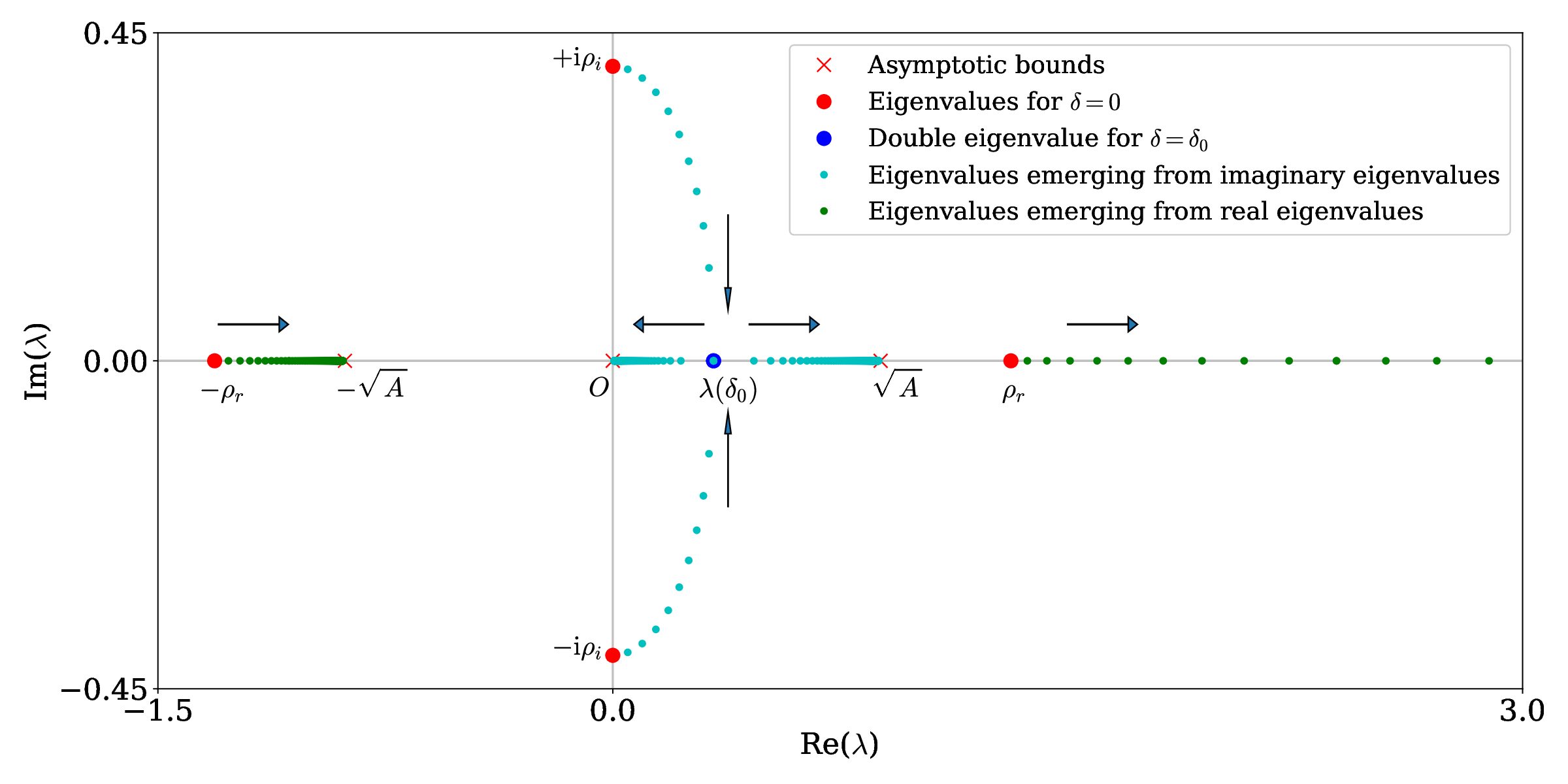}
\caption{Location of the eigenvalues of $\tilde{\bA}(\delta)$ in the complex plane for $s=1.2$ and various values of $\delta\geq 0$. The arrows indicate the direction of motion of the eigenvalues as $\delta$ increases from $0$.}
\label{fig:eigenvalues}
\end{figure}

By Proposition \ref{prop:eigenvalues}, for $0<\delta<\delta_0$ the leading eigenvalues consist of a pair of complex conjugate roots of $P_2(\lambda)$ with positive real part, together with one negative eigenvalue. As $\delta$ increases to $\delta_0$, the complex conjugate pair coalesces into a double positive real eigenvalue, which splits into two distinct positive real eigenvalues for $\delta>\delta_0$. Consequently, $\bU_0$ is a saddle-focus equilibrium for $0 < \delta < \delta_0$ and a saddle equilibrium for $\delta \geq \delta_0$ according to the definitions given in \cite{Kuz2024}.

Figure \ref{fig:eigenvalues} shows the eigenvalue locations for various values of $\delta$ with $s = 1.2$. Arrows indicate the trajectories of the eigenvalues in the complex plane as $\delta$ increases. The eigenvalues were computed numerically in Python, confirming that the equilibrium $\bU_0$ transitions from a saddle-focus for $0 < \delta < \delta_0$ to a saddle for $\delta \geq \delta_0$. This transition reflects a qualitative change in the shape of the heteroclinic orbit in the four-dimensional phase space.

It is worth noting that if $W_u(\bU)$ and $W_s(\bU)$ denote, respectively, the unstable and stable manifolds of an equilibrium $\bU$, then the Morse indices of the two equilibria in our system are $\ind(\bU_0) = \dim(W_u(\bU_0)) = 3$ and $\ind(\bo) = \dim(W_u(\bo)) = 2$. Therefore, a heteroclinic orbit can lie within the transverse intersections $W_u(\bU_0) \cap W_s(\bo)$, whose dimension is $\mathrm{ind}(\bU_0) - \mathrm{ind}(\bo) = 1$, cf. \cite{HS2010} and the references therein.

\subsection{Discussion on the properties of the solutions}\label{sec:discussion}

In the sequel, we assume that $\delta>0$. We consider the functional $H:\mathbb{R}^4\to \mathbb{R}$ such that
\begin{equation}\label{eq:hamiltonian}
    H(\bU)=Q(\bu')+f(\bu)=-su'\eta'+f(\bu)\ ,
\end{equation}
where $\bu=(u,\eta)$ and $\bu'=(u',\eta')$. 

Note that 
$$
\frac{d}{d\xi} f(\bu(\xi))=\nabla f(\bu)\cdot \begin{pmatrix} u'\\ \eta'\end{pmatrix}\ .
$$
Taking the inner product of (\ref{eq:sys2}) with the vector $(u',\eta')^T$ yields
\begin{equation}\label{eq:derivprod}
s(u'\eta')'=\frac{d}{d\xi}f(\bu)+\delta u'\eta'\ ,
\end{equation}
which can be written in the form
\begin{equation}\label{eq:potender}
\frac{d}{d\xi}H(\bU)=-\delta u' \eta' \ .    
\end{equation}

It is straightforward to verify using (\ref{eq:potender}) that for any solution of (\ref{eq:sys4d}) and for any $\xi\in\mathbb{R}$: 
\begin{equation}\label{eq:firstderiv}
    \frac{d}{d\xi}H(\bU)-\frac{\delta}{s}H(\bU)=-\frac{\delta}{s}f(\bu)\ .
\end{equation} 
Multiplying (\ref{eq:firstderiv}) with $e^{-\frac{\delta}{s}\xi}$, we obtain
\begin{equation}\label{eq:motion}
\frac{d}{d\xi}\left(e^{-\frac{\delta}{s}\xi}H(\bU)\right)=-\frac{\delta}{s}e^{-\frac{\delta}{s}\xi}f(\bu)\ .  
\end{equation}
If we denote 
\begin{equation}\label{eq:distance}
V(\xi)=e^{-\frac{\delta}{s}\xi}H(\bU(\xi))\ ,
\end{equation}
then we have
\begin{equation}\label{eq:velocity}
V'(\xi)=-\frac{\delta}{s}e^{-\frac{\delta}{s}\xi}f(\bu)\ ,
\end{equation}
and
\begin{equation}\label{eq:acceleration}
V''(\xi)=\frac{\delta}{s}e^{-\frac{\delta}{s}\xi}[\frac{\delta}{s}f(\bu)-\frac{d}{d\xi}f(\bu)]\ .
\end{equation}

In this work we focus on heteroclinic orbits connecting $\bU_0$ with the origin, where the domain $\mathcal{D}=\mathring{\mathcal{C}}\cap\{u>0,\eta>0\}$ seems to play a significant role.

\begin{remark}\label{rmk:exit}
If $\xi_0\in\mathbb{R}$ is a point at which the orbit $\bu$ crosses $\mathcal{C}$ from its interior, then, $V'(\xi_0)=0$, and for $\xi<\xi_0$, $V'(\xi)<0$ while for $\xi>\xi_0$, $V'(\xi)>0$. Thus, $V(\xi_0)$ is a local minimum, which implies that $V''(\xi_0)\geq 0$. Since $f(\bu(\xi_0))=0$, the definition (\ref{eq:acceleration}) implies that $\frac{d}{d\xi}f(\bu)|_{\xi=\xi_0}=\nabla f(\bu)\cdot \bu'|_{\xi=\xi_0}\leq 0$, as expected. Similarly, if $\xi_1$ is a point at which the orbit $\bu$ crosses $\mathcal{C}$ from its exterior, we have that $\frac{d}{d\xi}f(\bu)|_{\xi=\xi_1}=\nabla f(\bu)\cdot \bu'|_{\xi=\xi_1}\geq 0$. Crossings that are tangential to $\mathcal{C}$ cannot be excluded by this argument.
\end{remark}

System (\ref{eq:sys2}) also possesses first integrals: 
Solving (\ref{eq:derivprod}) for $u'\eta'$ (using integrating factors and the boundary conditions (\ref{eq:bc1})) in $(\xi,\infty)$, we obtain
\begin{equation}\label{eq:negative}
H(\bU)=\frac{\delta}{s} e^{\frac{\delta}{s}\xi}\int_\xi^\infty e^{-\frac{\delta}{s}\tau}f(\bu)~d\tau\ ,
\end{equation}
for all $\xi\in\mathbb{R}$. This, in turn, gives
$$V(\xi)=\frac{\delta}{s}\int_\xi^\infty e^{-\frac{\delta}{s}\tau} f(\bu)~d\tau \ ,$$
for all $\xi \in \mathbb{R}$. It follows that $V(\xi)$ (and $H(\bU)$) can vanish only if the orbit exits the domain $\mathcal{D}$, that is, if there exist intervals of $\xi$ where the function $f(\bu)$ becomes negative.

Integration of (\ref{eq:potender}) using the boundary conditions (\ref{eq:bc2}) leads to the first integral
\begin{equation}\label{eq:potender2}
H(\bU)-f(\bu_0)+\delta \int_{-\infty}^\xi u'\eta'~d\tau=0\ ,    
\end{equation}
which shows after taking the limit $\xi\to\infty$ that
\begin{equation}\label{eq:intinf}
\delta\int_{-\infty}^{+\infty}u'\eta'=f(\bu_0)>0\ ,
\end{equation}
which also yields from (\ref{eq:potender2})
\begin{equation}\label{eq:potender3}
H(\bU)=\delta \int_\xi^{+\infty} u'\eta'~d\tau\ .   
\end{equation}

\begin{remark}\label{rem:contra2}
The existence of homoclinic orbits to the origin inside the region $\mathcal{D}$ can be easily excluded. If there was such orbit, then $\displaystyle{\lim_{\xi\to-\infty}H(\bU)=0}$. Furthermore, there will be a $\rho>0$ such that $u'>0$ and $\eta'>0$ for all $\xi\in (-\infty,\rho)$, i.e. as the orbit progresses from the origin inside $\mathcal{D}$. On the other hand, (\ref{eq:potender}) yields that $\frac{d}{d\xi}H(\bU)<0$ in $(-\infty,\rho)$, and thus $H(\bU(\xi))<0$ for $\xi\in (-\infty,\rho)$, which contradicts with (\ref{eq:negative}) that ensures that $H(\bU)>0$ for all $\xi\in\mathbb{R}$.
\end{remark}

A direct consequence of the previous analysis is that any solution $\bU=(u,\eta,v,h)$ of the system (\ref{eq:sys4d}) with boundary conditions (\ref{eq:bc1})--(\ref{eq:bc2}) and $\bu\in \mathcal{D} =\mathring{\mathcal{C}}\cap\{u>0,\eta>0\}$ satisfies $0=H(\bo)<H(\bU)\leq  H(\bU_0)=f(\bu_0)$
for all $\xi\in\mathbb{R}$, and
\begin{equation}\label{eq:maximprinc}
\delta \int_{-\infty}^{+\infty} u'\eta'~d\tau=f(\bu_0) > 0\quad \text{and}\quad \int_{-\infty}^\xi u'\eta'~d\tau\geq 0,\quad \int_\xi^{+\infty} u'\eta'~d\tau> 0 \quad \text{for all $\xi\in\mathbb{R}$}\ ,
\end{equation}

To see this observe that if $\bu\in\mathcal{D}$ for all $\xi\in\mathbb{R}$, then (\ref{eq:negative}) implies that $H(\bU)>0$ for all $\xi\in\mathbb{R}$. Moreover, since $f(\bu)\leq f(\bu_0)$ for $\bu\in\mathcal{D}$, we have that the right-hand side of (\ref{eq:negative}) is 
$$\frac{\delta}{s} e^{\frac{\delta}{s}\xi}\int_\xi^\infty e^{-\frac{\delta}{s}\tau}f(\bu)~d\tau \leq f(\bu_0)\ ,$$
concluding that $0<H(\bU)\leq f(\bu_0)$. The second inequality in  (\ref{eq:maximprinc}) follows from (\ref{eq:potender2}) since $H(\bU)-f(\bu_0)\leq 0$.

On the other hand, heteroclinic orbits with $H(\bU)>0$ can only lie within $\mathcal{D}$, along the arc $P_1P_2$, and in the region bounded by the arcs $P_1P_2$, $P_2P_3$, and $P_1P_3$. Indeed, if $H(\bU)>0$, then $Q(\bu')>-f(\bu)$. Consequently, if an orbit exits $\mathcal{D}$, it must satisfy $Q(\bu')>0$, which prevents it from returning to $\mathcal{C}$ due to the convexity of the latter, except if the orbit exits $\mathcal{C}$ from the arc $P_1P_2$ but not $P_1P_3$ or $P_2P_3$.

It is noted that we cannot infer that critical points of solutions of (\ref{eq:sys2}) occur simultaneously. This was unexpected since for other systems, such as the dissipative classical Boussinesq--Peregrine system \cite{BMT2022}, the modified Serre--Green--Naghdi system \cite{BMT2024}, and the KdV--Burgers equation \cite{BS1985}, the critical points of the solutions $u$, $\eta$ occur simultaneously. Even if the maxima and minima of $u$ and $\eta$ do not occur simultaneously, solutions corresponding to undular bores are expected to have extrema of $u$ and $\eta$ that are very close to each other and of the same type. In other words, it is expected that the intervals between the roots of $u'$ and $\eta'$ where $u'\eta' < 0$ are small compared to the intervals where $u'\eta' > 0$. 

\begin{figure}[t]
\centering
\includegraphics[width=0.8
\columnwidth]{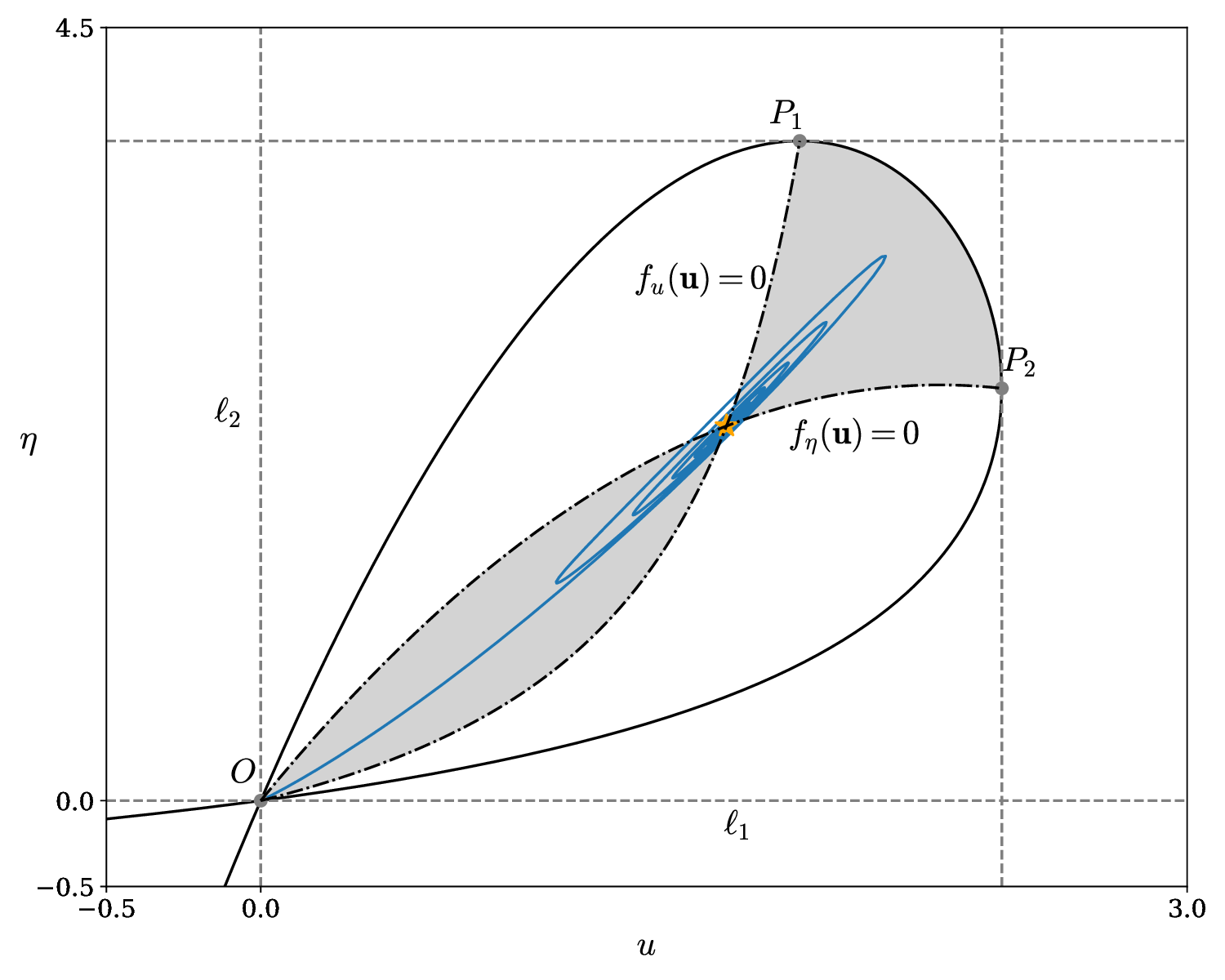}
\caption{Anticipated locations of the minima and maxima of $u$ and $\eta$ in $\mathcal{D}$ in physical space. The shaded region represents the set $\{\bu\in \mathcal{D}~:~f_u(\bu)f_\eta(\bu)>0\}$. In the shaded portion above $\bu_0$, both satisfy $f_u(\bu)<0$ and $f_\eta(\bu)<0$, whereas in the portion below $\bu_0$ we have $f_u(\bu)>0$ and $f_\eta(\bu)>0$. Here, $s=2.2$ and $\delta=0.5$.}
\label{fig:minmax}
\end{figure}

Consider a point $\xi_0$ where $u(\xi_0)$ is a local maximum of $u$, i.e. $u'(\xi_0)=0$ and $u''(\xi_0)<0$. Then, from the second equation of (\ref{eq:sys1}) for $\xi=\xi_0$ we have that $\left. f_\eta(\bu)\right|_{\xi=\xi_0}=su(\xi_0)-\eta(\xi_0)-\tfrac{1}{2}u^2(\xi_0)<0$. If $\xi_1$ is a point where the solution $u$ achieves a minimum, then 
$\left. f_\eta(\bu)\right|_{\xi=\xi_1}=su(\xi_1)-\eta(\xi_1)-\tfrac{1}{2}u^2(\xi_1)>0$. Similarly, if at $\xi=\zeta_0$ the solution $\eta(\zeta_0)$ is a local maximum of $\eta$, then $\left. f_u(\bu)\right|_{\xi=\zeta_0}=s\eta(\zeta_0)-u(\zeta_0)-\eta(\zeta_0)u(\zeta_0)<0$, and if $\xi=\zeta_1$ is a point where $\eta(\zeta_1)$ is a local minimum, then $\left. f_u(\bu)\right|_{\xi=\zeta_1}=s\eta(\zeta_1)-u(\zeta_1)-\eta(\zeta_1)u(\zeta_1)>0$.
For physical solutions that lie entirely inside $\mathcal{D}$, as mentioned earlier, the points where the maxima of $\eta$ and $u$ occur are expected to be very close to each other. The same holds for their minima. Consequently, the region inside $\mathcal{D}$ where extrema of $u$ and $\eta$ can occur is expected to correspond to the set of points satisfying where $f_u(\bu) \cdot f_\eta(\bu) > 0$. This region is shown in Figure~\ref{fig:minmax}. We verified this claim numerically for all the cases we tested. However, due to finite-precision arithmetic, it is not possible to confirm it for very small values of $\delta$ and $\xi$.

In addition to the above, Equation (\ref{eq:potender}) implies that if $\xi_0 \in \mathbb{R}$ is a point such that $u'(\xi_0)\eta'(\xi_0) = 0$, then $\frac{d}{d\xi}H(\bU) = 0$ at $\xi = \xi_0$. Thus, the critical points of $H$ coincide with those of $\eta$ and $u$. Moreover, from the definition (\ref{eq:hamiltonian}) of $H$, we see that at the critical points of $u$ and $\eta$, we have $H(\bU) = f(\bu)$. Therefore, a local extremum of $H$ is always less than $f(\bu_0)$, provided that $\bu$ lies inside $\mathcal{D}$, or even in a broader region where $f(\bu) < f(\bu_0)$. Consequently, at least initially, $H(\bU)$ will be decreasing and $u'\eta' > 0$. Similarly, if $\bu \in \mathcal{D}$, then the local extrema of $H$ will be positive.

\section{Dissipative Bona--Smith-type systems}\label{sec:bonasmith}

Here we consider the Bona--Smith system, 
\begin{equation}\label{eq:BS1}
\begin{aligned}
    &\eta_t+u_x+(u\eta)_x-\eta_{xxt}=0\ ,\\
    &u_t+\eta_x+uu_x-\gamma \eta_{xxx}-u_{xxt}=\delta u_{xx}\ ,
\end{aligned}
\end{equation}
where $\gamma\in (0,1]$. Similarly to the BBM--BBM system, we search for traveling wave solutions in the form $u(x,t)=u(x-st)$, $\eta(x,t)=\eta(x-st)$, $s>1$. Substitution into (\ref{eq:BS1}) leads to the system
\begin{equation}\label{eq:bsys1}
    \begin{aligned}
        &-s\eta''+s\eta-u-\eta u=0\, \\
        &-su''+\gamma \eta''+su-\eta-\tfrac{1}{2}u^2=-\delta u'\ ,
    \end{aligned}
\end{equation}
which is of the form (\ref{eq:sys2}) with
$$\bS=\begin{pmatrix}
    0 & -s \\ 
    -s & \gamma
\end{pmatrix}, \quad \bK=\begin{pmatrix}
    0& 0 \\
    -1 & 0
\end{pmatrix}\quad \text{and} \quad f(\bu)=s~u\eta-\tfrac{1}{2}\left(\eta^2+u^2+u^2\eta \right)\ .$$

The equilibrium solutions $\bo$ and $\bu_0=(u_0, \eta_0)$ of system (\ref{eq:bsys1}) remain the same with the BBM--BBM system for all $\gamma\in [0,1]$. The Jacobian matrix of the corresponding first order system linearized around any point $(u,\eta,0,0)$ is
\begin{equation}\label{eq:bslin1}
\bJ(\delta;\bu)=\begin{pmatrix}
0 & 0 & \tfrac{1}{s} & 0\\
0 & 0 & 0 & \tfrac{1}{s}\\
s-u-\tfrac{\gamma}{s}(1+\eta) & -1+\tfrac{\gamma}{s}(s -u) & \tfrac{\delta}{s} & 0 \\
-(1+\eta) & s-u & 0 & 0
\end{pmatrix}\ .
\end{equation}
Thus, the linearization about the origin $\bo$ corresponds to the matrix $\bA(\delta)=\bJ(\delta;\bo)$ with characteristic polynomial 
\begin{equation}
P_1(\lambda)=\lambda^4-\frac{\delta}{s}\lambda^3-\left(2-\frac{\gamma}{s^2}\right)\lambda^2+\frac{\delta}{s}\lambda+\left(1-\frac{1}{s^2}\right)\ .
\end{equation}
Again, using Descartes' rule of signs, we conclude that it has two positive and two negative roots. It is remarkable that in the case $\gamma=1$ two of the eigenvalues are the $\pm 1$ for all $s>1$ since $P_1(\pm 1)=(\gamma-1)/s^2$.

The unperturbed matrix $\bA(0)$ for $\delta=0$ has four real eigenvalues (similar to the BBM--BBM case), namely,
$$\begin{aligned}&\rho_1=-\sqrt{1+\frac{\Gamma-\gamma}{2s^2}}\leq-1,\quad \rho^{-}=-\sqrt{1-\frac{\Gamma+\gamma}{2s^2}}>-1\ ,\\
& \rho^+=\sqrt{1-\frac{\Gamma+\gamma}{2s^2}}<1,\quad \rho_2=\sqrt{1+\frac{\Gamma-\gamma}{2s^2}}\geq 1\ ,\end{aligned}$$
where $\Gamma=\sqrt{\gamma^2+4s^2(1-\gamma)}\geq\gamma$, with equality only for $\gamma=1$, in which case $\rho_1=-1$ and $\rho_2=1$. It is easy to verify that for $\gamma=0$ these eigenvalues coincide with the corresponding eigenvalues of the BBM--BBM system. 

The eigenvalues of the perturbed matrix $\bA(\delta)=\bJ(\delta;\bo)$ as $\delta\to 0$ behave like
$$\lambda_1\sim\rho_1+\frac{\delta}{4s}\frac{\Gamma-\gamma}{\Gamma},\quad \lambda^-\sim \rho^{-}+\frac{\delta}{4s}\frac{\Gamma+\gamma}{\Gamma}, \quad \lambda^{+}\sim \rho^{+}+\frac{\delta}{4s}\frac{\Gamma+\gamma}{\Gamma}, \quad \lambda_2\sim \rho_2+\frac{\delta}{4s}\frac{\Gamma-\gamma}{\Gamma}\ ,$$
while as $\delta\to\infty$ behave like
$$\lambda_1\sim -1+\frac{\gamma-1}{2\delta s},\quad \lambda^{-}\sim -\frac{s^2-1}{\delta s}, \quad \lambda^{+}\sim 1+\frac{\gamma-1}{2\delta s},\quad \lambda_2\sim \frac{\delta}{s}\ .$$

The corresponding linearization matrix $\tilde{\bA}(\delta)=\bJ(\delta;\bu_0)$ about the equilibrium $\bu_0$  has the following characteristic polynomial
\begin{equation}
P_2(\lambda)=\lambda^4-\frac{\delta}{s}\lambda^3-2\left(\frac{s-u_0}{s}-\frac{\gamma}{2}\frac{1+\eta_0}{s^2} \right)\lambda^2+\frac{\delta}{s}\frac{s-u_0}{s}\lambda-\frac{1+\eta_0-(s-u_0)^2}{s^2}\ .
\end{equation}
The coefficient $\tilde{A}$ of the quadratic term is
\begin{equation}\label{eq:Aformula2}
\tilde{A}=A-B=\frac{s-u_0}{s}-\frac{\gamma}{2}\frac{1+\eta_0}{s^2} =-\frac{\gamma(s+\sqrt{s^2+8})+4(s-\sqrt{s^2+8})}{8 s}\ ,
\end{equation} where
$A$ is given by (\ref{eq:Aformula1}) and 
$$B=\frac{\gamma}{2}\frac{1+\eta_0}{s^2}\ .$$
Observe that $\tilde{A}>0$ if and only if $\gamma(s+\sqrt{s^2+8})+4(s-\sqrt{s^2+8})<0$, i.e. if and only if $1<s<\frac{4-\gamma}{\sqrt{2\gamma}}$. For larger values of $s$ we have $\tilde{A}\leq 0$. For $\gamma=0$ we have $\tilde{A}=A>0$ for all $s>1$.

Applying the Descartes' rule of signs, (similarly with the BBM--BBM system), we deduce that $P_2(\lambda)$ has 1 positive and 1 negative eigenvalues, and depending on the value of $\delta$, we can have additionally either 2 complex conjugate (with positive real part) or 2 real positive eigenvalues. 

When $\delta=0$, the unperturbed system has two real $\rho_r^\pm=\pm\rho_r$ and two imaginary eigenvalues $\rho_i^\pm=\pm\I\rho_i$, where
$$
\rho_i=\sqrt{\frac{u_0-s}{s}+\gamma\frac{1+\eta_0}{2s^2}+\sqrt{1+\eta_0}\frac{\Theta}{2s^2}},\quad 
    \rho_r=\sqrt{\frac{s-u_0}{s}-\gamma\frac{1+\eta_0}{2s^2}+\sqrt{1+\eta_0}\frac{\Theta}{2s^2}}\ ,
$$
where $\Theta=\sqrt{\gamma^2(1+\eta_0)+4s^2(1-\gamma)+4\gamma su_0}$.

Denoting the real eigenvalues by $\lambda_r^{\pm}$ emerging from the real eigenvalues $\rho_r^\pm=\pm\rho_r$ of the unperturbed system ($\delta=0)$ and $\lambda_i^{\pm}$ those that emerge from the imaginary eigenvalues $\rho_i^\pm$ of the unperturbed system, we can compute the asymptotic limits as $\delta\to 0$
$$\lambda_r^{\pm}\sim \frac{\delta}{4s}\frac{\Theta-\gamma\sqrt{1+\eta_0}}{\Theta}\pm\rho_r+O(\delta^2) ,\quad \lambda_i^{\pm}\sim\frac{\delta}{4s}\frac{\gamma\sqrt{1+\eta_0}+\Theta}{\Theta}\pm\I\rho_i+O(\delta^2)\ , $$
and as $\delta\to\infty$
$$\begin{aligned}
&\lambda_i^{-}\sim \frac{1+\eta_0-(s-u_0)^2}{\delta(s-u_0)},\quad \lambda_i^+\sim\sqrt{\frac{s-u_0}{s}}-\frac{(1+\eta_0)[(1-\gamma)s+\gamma u_0]}{2 \delta s(s-u_0)}\ ,\\
&\lambda_r^{-}\sim  -\sqrt{\frac{s-u_0}{s}}-\frac{(1+\eta_0)[(1-\gamma)s+\gamma u_0]}{2\delta s(s-u_0)},\quad \lambda_r^+\sim \frac{\delta}{s}\ .
\end{aligned} $$
It is worth mentioning that the previous relations for $\gamma=0$ coincide with those derived in the case of the BBM--BBM equation.

Note that (\ref{eq:eigsderiv}), which was derived in the case of the BBM--BBM system, is the same in the case of Bona--Smith systems, since the coefficients of $\lambda^3$, $\lambda$ and $\lambda^0$ of $P_2(\lambda)$ are not affected by $\gamma$. Moreover, $P_2(\pm\sqrt{A})=-\frac{(1+\eta_0)[(1-\gamma)s+\gamma u_0]}{s^3}\neq 0$. However, the coefficient of $\lambda^2$ is now $-2\tilde{A}$ instead of $-2A$, and $\tilde{A}$ can be negative. For this reason, the arguments of Appendix \ref{app:linearization} that rely on the inflection points of $P_2(\lambda)$ and on the specific form of the coefficient of $\lambda^2$ do not carry over verbatim. Formula (\ref{eq:Xofm}) remains valid, and thus $Re(\lambda_i^\pm)>0$ for all $\delta>0$. The remaining steps of the analysis of Appendix \ref{app:linearization} (the monotonicity of $m(\delta)$, the collision of the complex eigenvalues at a double positive eigenvalue and the monotonicity of their imaginary parts) were verified numerically for $\gamma\in\{0.25,0.5,0.75,1\}$ and $1<s\leq 60$, and the details are omitted.

To determine the critical bifurcation value $\delta_0$, we compute the last remainder in Euclid's algorithm for the greatest common divisor of $P_2(\lambda)$ and $P_2'(\lambda)$. This remainder is found to be
$$\frac{\left(16 s^2 (A-B)+3 \delta^2\right)^2 \tilde{S}(\delta)}{64 \left(2 \delta^2 s^3 \left(7 A^2-18 A B+2 B^2+3 q\right)+32 s^5 (A-B) \left((A-B)^2+q\right)+3 A \delta^4 s\right)^2}\ ,$$
where $A$ and $q$ as before, and
$$\begin{aligned}
\tilde{S}(\delta)=& 4 A^3 \delta^6+  \left(A^2 \left(13 A^2-44 A B+4 B^2\right)+6 A q (6 B-5 A)-27 q^2\right)s^2 \delta^4\\
&~+32  \left(-q \left(2 A^2-12 A B+B^2\right) (A-B)+A^2 (A-B)^3-3 q^2 (A-3 B)\right)s^4 \delta^2 \\
&~ -256 q  \left((A-B)^2+q\right)^2 s^6\ .
\end{aligned}
$$
We conclude that the complex eigenvalues collapse into a double positive eigenvalue when $\tilde{S}(\delta) = 0$. By Descartes' rule of signs, the corresponding cubic polynomial in the variable $x = \delta^2 / s^2$ has at least one positive real root. To verify that this root is unique, we computed symbolically (using Mathematica\textsuperscript{\tiny\textregistered}) the discriminant of the cubic polynomial $\tilde{S}(x)=ax^3+bx^2+cx+d$ given by
$$\Delta=18abcd +b^2c^2-4ac^3-4b^3d-27a^2d^2\ ,$$ and confirmed that $\Delta < 0$ for all $s > 1$ and $\gamma \in [0,1]$. The resulting expressions, however, are not particularly elegant, and we omit them here.

In addition to these symbolic calculations, we numerically studied the roots of the cubic polynomial $\tilde{S}(x)$ for several values of $\gamma$ and $s > 1$. In all cases, we found exactly one positive real root and two complex conjugate roots. We denote the positive real root of $\tilde{S}(\delta) = 0$ (the critical bifurcation value) by $\delta_0(\gamma)$, to emphasize its dependence on $\gamma$. For $\gamma = 0$, the value $\delta_0(0)$ coincides with that obtained in Section \ref{sec:dynamics}.

Finally, we note that the energy functional $H$ in this case takes the form
\begin{equation}
H(\bU)=-su'\eta' + \frac{\gamma}{2}(\eta')^2 + f(\bu)\ .
\end{equation}
Because the definition of $H$ differs between the BBM--BBM and Bona--Smith systems, the behavior of critical points of $u$ and $\eta$ is also different. In the Bona--Smith case, a critical point of $\eta$ coincides with a critical point of $H$, and we have $H(\bU)<f(\bu_0)$. By contrast, at a critical point of $u$ one obtains $H(\bU)=\frac{\gamma}{2}(\eta')^2+f(\bu)$, which is different from the BBM--BBM situation. Consequently, the regions in which maxima and minima of the solutions may occur are not the same in the two systems, with the maxima and minima of the traveling waves of the Bona--Smith system to be bounded by the curves $f_u(\bu)=0$ and $\tfrac{\gamma}{s}f_u(\bu)+f_\eta(\bu)=0$. However, we verified experimentally that for sufficiently large values of $\xi$, the maxima and minima were in fact confined to the same regions. Moreover, we confirmed experimentally that $H(\bU)<f(\bu_0)$ for all computed $\xi$, and even though $H(\bU)$ in some cases approached $f(\bu_0)$ in the experiments of Section \ref{sec:comparison}, the distance $|H(\bU)-f(\bu_0)|$ consistently remained greater than $10^{-3}$.

\begin{figure}[t]
\centering
\includegraphics[width=0.7\columnwidth]{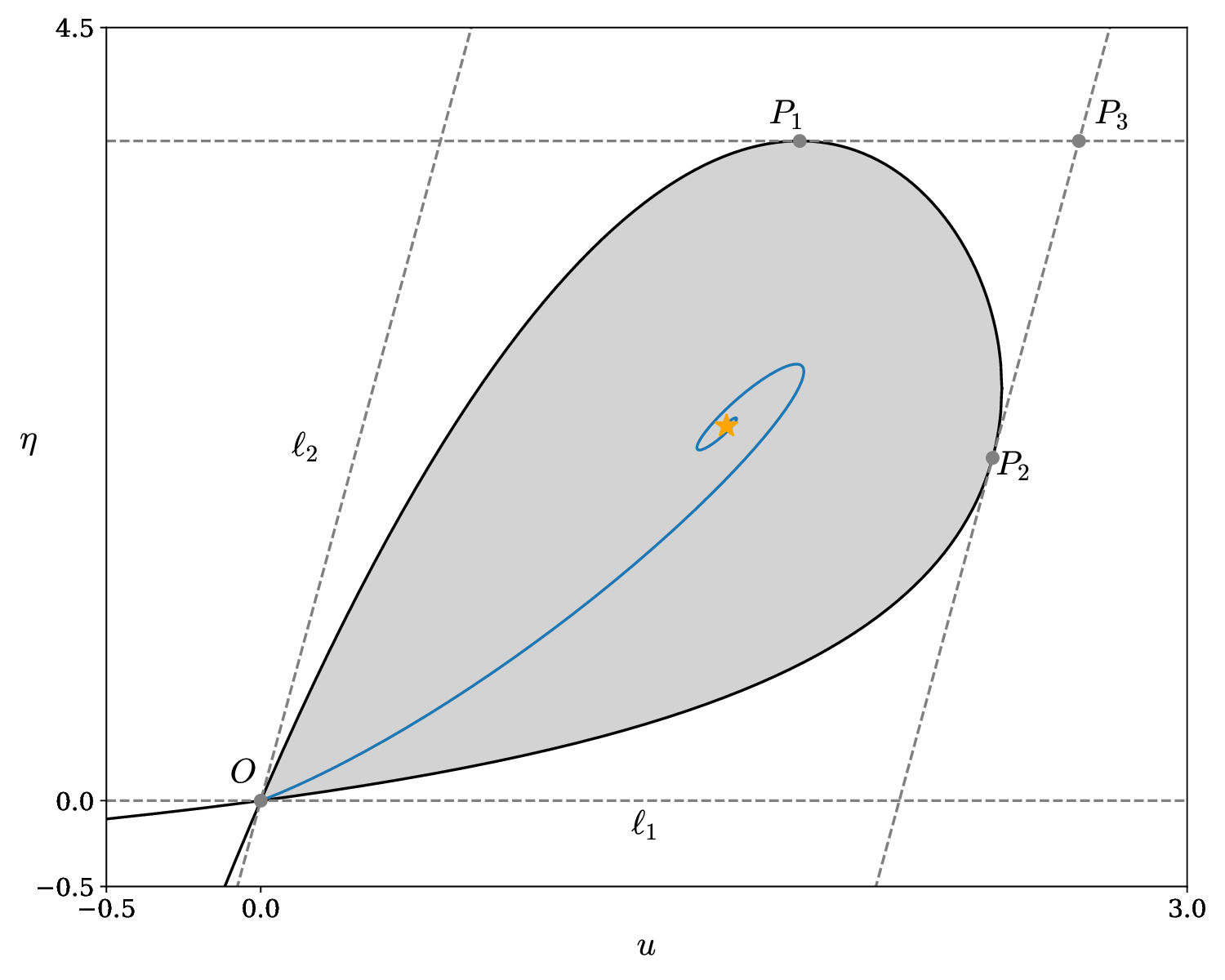}
\caption{A projection of a heteroclinic orbit of the Bona--Smith system with $\gamma=2/3$ connecting the critical points $\bo$ and $\bU_0$ (depicted by a star) on physical space when $s=2.2$ and $\delta=1.8$. The projected orbit $\mathcal{R}$ lies entirely inside the interior of the shaded domain $\mathcal{D}$.}
\label{fig:domain2}
\end{figure}

Finally, we note that for the Bona--Smith system the line $\ell_2$ is given by $\eta=\tfrac{2s}{\gamma}u$. Together with the line $\ell_1$ (the $u$-axis), it bounds the region where $Q(\bu)<0$, within which the curve $\mathcal{C}$ lies in the interior; see Figure~\ref{fig:domain2}. Since in this case $-su'\eta'=H(\bU)-f(\bu)-\frac{\gamma}{2}(\eta')^2$, the relations (\ref{eq:firstderiv})--(\ref{eq:negative}) of Section \ref{sec:discussion} hold true for the Bona--Smith system after replacing $f(\bu)$ with $f(\bu)+\frac{\gamma}{2}(\eta')^2$. In particular,
$$H(\bU)=\frac{\delta}{s} e^{\frac{\delta}{s}\xi}\int_\xi^\infty e^{-\frac{\delta}{s}\tau}\left(f(\bu)+\frac{\gamma}{2}(\eta')^2\right)d\tau\ ,$$
and thus $H(\bU)>0$ for any heteroclinic orbit that lies in $\mathcal{D}$. With this modification, the rest of the conclusions mentioned in Section \ref{sec:discussion} hold true also for the Bona--Smith system.

\section{Numerical computation of the heteroclinic orbits}\label{sec:computho}

In this section we present a numerical method for the computation of heteroclinic orbits that resemble oscillatory and regularized shock waves of the aforementioned Boussinesq-type systems. For the sake of simplicity, we present the method only for the BBM--BBM system, however, similar procedure was used for the Bona--Smith as well as the classical Boussinesq--Peregrine system.

Via \eqref{eq:sys1} we may write
\begin{equation*}
\eta = -su''+su-\frac{1}{2}u^2+\delta u' \ ,
\end{equation*}
and via subsequent substitution we find
\begin{equation}\label{eq:sysu}
s^2u^{(4)}-\delta s u^{(3)}-2s(s-u)u'' + (\delta(s-u)+su')u'+(s^2-1-\frac{3}{2}su+\frac{1}{2}u^2)u = 0 \ .
\end{equation}
We seek solutions for which $u\to u_0$ and $u',u'',u^{(3)}\to 0$ as $\xi\to-\infty$ and $u,u',u'',u^{(3)}\to 0$ as $\xi\to+\infty$.
We consider the asymptotic solution at both extremes, a potentially useful mechanism for matching with numerical solutions on a finite interval.

\begin{remark} 
For the more general family of Bona--Smith systems \eqref{eq:bsys1} is similarly reduced to a fourth order problem in $u$ by first combining the two equations to eliminate $\eta''$ and rearranging the result to obtain an expression for $\eta$ in terms of $u,u',u''$. 
Substituting the resulting expression into either of the original two equations, followed by some algebraic manipulations, leads to a (non-linear) fourth order ODE for $u$, specifically
\begin{align*}
0=&-s\left(s^2u^{(4)}+(\gamma-s^2)u''+su'^2+suu''-\delta su^{(3)}\right)(\gamma(s-u)-s)^2 \\
&-s\gamma u''\left(s^2u''+(\gamma-s^2)u+\tfrac{1}{2}su^2-\delta s u'\right)(\gamma(s-u)-s) \\
&-2s\gamma u'\left(s^2u^{(3)}+(\gamma-s^2)u'+suu'-\delta su^{(2)}\right)(\gamma(s-u)-s) \\
&-2s\gamma^2 u'^2 \left(s^2u''+(\gamma-s^2)u+\tfrac{1}{2}su^2-\delta s u'\right) \\
&+(s-u)\left(s^2u''+(\gamma-s^2)u+\tfrac{1}{2}su^2-\delta s u'\right)(\gamma(s-u)-s)^2
-u(\gamma(s-u)-s)^3 \,.
\end{align*}
This is then solved analogously to what is described here for the BBM--BBM system.
\end{remark}

Consider first the case $\xi\to+\infty$.
Since each of $u,u',u'',u^{(3)}$ decays to zero, then the leading order component of the asymptotic behavior solves the linearized problem
\begin{equation*}
s^2u^{(4)}-\delta s u^{(3)}-2s^2 u'' + \delta s u'+(s^2-1)u = 0 \,, \quad \text{as}\quad  \xi\to+\infty\ .
\end{equation*}
The characteristic equation for this linear constant coefficient differential equation is equivalent to $s^2P_1(\lambda)$, see \eqref{eq:P1poly}.
Thus, we expect exactly two of the eigenvalues to be negative and exactly two to be positive.
The negative eigenvalues are the ones which result in the desired decay of the solution as $\xi\to\infty$, so denoting them by $\lambda_{+1},\lambda_{+2}$ the asymptotic solution must have the form
\begin{equation}
u \sim u_{+\infty} := c_{+1}e^{\lambda_{+1}\xi}+c_{+2}e^{\lambda_{+2}\xi}, \quad \xi\to\infty \ .
\end{equation}

Analogously, we consider the case $\xi\to-\infty$.
It is convenient here to introduce $u=u_0+\tilde{u}$ such that $\tilde{u},\tilde{u}',\tilde{u}'',\tilde{u}^{(3)}\to 0$ as $\xi\to-\infty$ and the leading order asymptotic behavior solves the linearized problem 
\begin{equation*}
s^2\tilde{u}^{(4)}-\delta s \tilde{u}^{(3)}-2s(s-u_0) \tilde{u}'' + \delta (s-u_0) \tilde{u}'+(u_0-\frac{3}{2}s)u_0\tilde{u} = 0 \,, \quad\text{as} \quad \xi\to-\infty\ .
\end{equation*}
The characteristic equation for this linear constant coefficient differential equation is equivalent to $s^2P_2(\lambda)$, see \eqref{eq:P2poly}.
Thus, with $\delta>0$ we expect exactly three eigenvalues to have positive real part (two of the three being complex for $0<\delta<\delta_0$) while the fourth eigenvalue is always negative real.
Those eigenvalues with positive real part are the ones which have the desired decay as $\xi\to-\infty$, so denoting them by $\lambda_{-1},\lambda_{-2},\lambda_{-3}$, when no two are equal the asymptotic solution must have the form
\begin{equation}
u \sim u_{-\infty} := u_0+c_{-1}e^{\lambda_{-1}\xi}+c_{-2}e^{\lambda_{-2}\xi}+c_{-3}e^{\lambda_{-3}\xi}, \quad \xi\to-\infty \ .
\end{equation}
When $\delta<\delta_0$ two of the eigenvalues are a complex conjugate pair, and as we seek real solution the corresponding coefficients are also complex conjugate.
In the case $\delta=\delta_0$, for which there is a double root of $P_2(\lambda)$ and without loss of generality we let $\lambda_{-2}=\lambda_{-3}$, the asymptotic solution has the form
\begin{equation*}
u \sim u_{-\infty} := u_0+c_{-1}e^{\lambda_{-1}\xi}+(c_{-2}+c_{-3}\xi)e^{\lambda_{-2}\xi}, \quad \xi\to-\infty \ .
\end{equation*}

\begin{remark}
In the case of the classical Boussinesq--Peregrine system the desired heteroclinic orbit is described via the solution of the second order problem 
$$-su(u-s)-u+\tfrac{1}{2}u^2(u-s)+ s(u-s)u''-\bar{\delta} (u-s)u'=0\,.$$
We once more utilize asymptotic solutions at $\xi=\pm\infty$ to solve this problem.
In this case the eigenvalues associated with the characteristic polynomial of the linearized problem in the limit $\xi\to+\infty$ always has exactly one negative (real) eigenvalue.
In the limit $\xi\to-\infty$ the eigenvalues associated with the linearized problem once more depend on the value of $\bar{\delta}$ relative to a specific $\bar{\delta}_0$, specifically there are: two distinct positive real eigenvalues when $\bar{\delta}>\bar{\delta}_0$; a repeated (double) positive real eigenvalue when $\bar{\delta}=\bar{\delta}_0$; and two complex eigenvalues with positive real part when $0<\bar{\delta}<\bar{\delta}_0$.
The asymptotic solutions are otherwise utilized in a similar way to what is described here for the BBM--BBM system.
\end{remark}

One difficulty with solving \eqref{eq:sysu} is that any translation of a given solution must also be a solution.
Among other things, this creates an issue in which an algebraic approximation of the system is not closed.
Noting that we expect $u'$ to possess a global minimum, i.e. a point of inflection where $u$ is decreasing steepest, we can make use of this to close the system and ``anchor'' the solution.
By way of example, by enforcing $u''(\xi_0)=0$ at some specific $\xi_0$ (nominally $\xi_0=0$) we may combine the two asymptotic solutions to construct an initial guess for any iterative method used to solve the nonlinear problem \eqref{eq:sysu}.
To illustrate, the two asymptotic solutions together have exactly five unknown coefficients. 
We can match them at the origin by enforcing the five conditions:
\begin{equation*}
u_{-\infty}(0)=u_{+\infty}(0) \,, \quad
u_{-\infty}^\prime(0)=u_{+\infty}^\prime(0) \,, \quad
u_{-\infty}^{\prime\prime}(0) = 0 \,, \quad
u_{+\infty}^{\prime\prime}(0) = 0 \,, \quad
u_{-\infty}^{(3)}(0)=u_{+\infty}^{(3)}(0) \,.
\end{equation*}
Having solved for the five unknown coefficients, we obtain an approximate solution to \eqref{eq:sysu} which has the anticipated asymptotic behavior at each end.
Through the course of subsequent iterations (e.g. formed via Newton's method), we can continue to enforce $u''(0)=0$ as a means of closing the system whilst keeping the interesting parts of the solution ``anchored'' near the origin.
Observe that in doing so we make no assumption about whether or not $\eta''(0)$ may be simultaneously zero.

Another challenge which arises is that the eigenvalues associated with the asymptotic behavior can have vastly different magnitude, especially in cases $\delta\ll \delta_0$ or $\delta\gg\delta_0$.
This can lead to an ill-conditioned system and some care is required.

Our implementation consists of deriving a Newton iteration associated with \eqref{eq:sysu} and then discretize the resulting iteration using centered 2nd order finite differences on a finite interval about the origin (noting the initial guess formed from the asymptotic approximations is useful for informing the choice of a suitable finite interval).
The corresponding algebraic system is then augmented to include the five (unknown) coefficients associated with the asymptotic solutions at each end.
The fourth order differential equation requires 4 boundary conditions and we require five additional pieces of information (one for each coefficient of the asymptotic solutions) to close the system.
We enforce that one-sided finite difference approximations of $u,u',u'',u^{(3)}$ match the appropriate asymptotic solution at each end point, thus constituting $4\times2=8$ pieces of information.
Some rescaling of parameters is performed to improve the conditioning of the equations associated with these matching conditions at the two ends.
For the final piece of information, we enforce that $u''(0)=0$ via the standard centered finite difference approximation of $u''$. 
Starting with the initial guess obtained from the asymptotic solutions, we achieve convergence reasonably quickly, typically less than $10^{-6}$ relative error (in the $2$-norm) in as few as four iterations (provided $\delta$ is not too small nor too large).

\subsection*{Non-linear iteration details}

Suppose that $u_k$ is some approximation to $u$ and we wish to find a correction $\Delta u$ such that $u_{k+1}=u_k+\Delta u$.
Substituting into \eqref{eq:sysu} we obtain
\begin{multline}\label{eq:newtoniteration}
s^2(\Delta u)^{(4)} -\delta s(\Delta u)^{(3)}-2s(s-u_k)(\Delta u)''  \\
+(\delta(s-u_k)+2su_k^\prime)(\Delta u)' +(s^2-1-3su_k+\frac{3}{2}u_k^2-\delta u_k^\prime+2su_k^{\prime\prime})(\Delta u) +O((\Delta u)^2) \\
= -s^2u_k^{(4)} +\delta s u_k^{(3)} +2s(s-u_k)u_k^{\prime\prime} -(\delta(s-u_k)+su_k^\prime)u_k^\prime -(s^2-1-\frac{3}{2}su_k+\frac{1}{2}u_k^2)u_k \,.
\end{multline}
Ignoring those $O((\Delta u)^2)$ terms, the solution to this equation provides an approximation to $\Delta u$ which constitutes a Newton method for solving \eqref{eq:sysu}.
To start off this iteration, we choose $u_1$ to be the approximation based on the two asymptotic solutions matched at $\xi=0$, as described above.

The nature of $u_1$ provides some insight into a suitable interval $[a,b]$, with $a<0<b$, for which the derivatives of $u$ have sufficiently decayed at the endpoints.
We then discretize this interval into equidistant points $\xi_i=a+(b-a)i/n$ for $i=0,1,2,\dots,n$ at which each of $u_k,u_{k+1},\Delta u$ will be estimated, along with their derivatives, using finite differences.
Specifically, in this way we discretize \eqref{eq:newtoniteration} at $\xi_2,\xi_3,\dots,\xi_{n-3},\xi_{n-2}$ using centered finite difference stencils.
For a typical fourth order boundary problem we would require four boundary conditions which, for example, could be implemented at the points $\xi_0,\xi_1,\xi_{n-1},\xi_n$.
However, for this problem we wish to ensure
\begin{align}\label{eq:asymptoticmatching}
u^{(j)}(a) &= u_{-\infty}^{(j)}(a) \,, \\
u^{(j)}(b) &= u_{+\infty}^{(j)}(b) \,, \notag
\end{align}
for each $j=0,1,2,3$. 
Recalling that the asymptotic solutions collectively have five unknown coefficients associated with them, these too need to be treated appropriately with respect to the Newton iteration, 
e.g. at $\xi=a$ we expand the matching conditions to obtain
\begin{align*}
u_k(a)+(\Delta u)(a) &= u_0+(c_{-3}+\Delta c_{-3})e^{\lambda_{-3}a}+(c_{-2}+\Delta c_{-2})e^{\lambda_{-2}a}+(c_{-1}+\Delta c_{-1})e^{\lambda_{-1}a} \,, \\
u_k^\prime(a)+(\Delta u)^{\prime}(a) &= (c_{-3}+\Delta c_{-3})\lambda_{-3}e^{\lambda_{-3}a}+(c_{-2}+\Delta c_{-2})\lambda_{-2}e^{\lambda_{-2}a}+(c_{-1}+\Delta c_{-1})\lambda_{-1}e^{\lambda_{-1}a} \,, \\
u_k^{\prime\prime}(a)+(\Delta u)^{\prime\prime}(a) &= (c_{-3}+\Delta c_{-3})\lambda_{-3}^2e^{\lambda_{-3}a}+(c_{-2}+\Delta c_{-2})\lambda_{-2}^2e^{\lambda_{-2}a}+(c_{-1}+\Delta c_{-1})\lambda_{-1}^2e^{\lambda_{-1}a} \,, \\
u_k^{(3)}(a)+(\Delta u)^{(3)}(a) &= (c_{-3}+\Delta c_{-3})\lambda_{-3}^3e^{\lambda_{-3}a}+(c_{-2}+\Delta c_{-2})\lambda_{-2}^3e^{\lambda_{-2}a}+(c_{-1}+\Delta c_{-1})\lambda_{-1}^3e^{\lambda_{-1}a} \,, 
\end{align*}
and similarly at $\xi=b$ (noting that here $u_0$ is the constant of \eqref{eq:equilibria}, not an iteration).
The derivatives associated with $u_k$ and $\Delta u$ at $\xi=a,b$ are each discretized using one-sided second order finite difference stencils.
Additionally, letting $i_0$ be an index which minimizes $|\xi_i|$ (i.e. ideally $\xi_{i_0}=0$, but is otherwise sufficiently close), then we also wish to ensure 
\begin{equation}\label{eq:anchoringcondition}
u_k^{\prime\prime}(\xi_{i_0})+(\Delta u)''(\xi_{i_0}) = 0 \,.
\end{equation}
Here the derivatives $u_k^{\prime\prime},(\Delta u)^{\prime\prime}$ are approximated using the standard centered finite difference stencil.
We then construct and solve the linear system
\begin{equation*}
\begin{bmatrix}
\text{row}_1  & \cdots \\
\text{row}_2  & \cdots \\
[M] & [0_{(n-3)\times 5}] \\
\text{row}_n  & \cdots \\
\text{row}_{n+1}  & \cdots \\
\text{row}_{n+2}  & \cdots \\
\text{row}_{n+3}  & \cdots \\
\text{row}_{n+4}  & \cdots \\
\text{row}_{n+5}  & \cdots \\
\text{row}_{n+6}  & \cdots
\end{bmatrix}
\begin{bmatrix}
(\Delta u)(\xi_0) \\
(\Delta u)(\xi_1) \\
\vdots \\
(\Delta u)(\xi_{n-1}) \\
(\Delta u)(\xi_n) \\
(\Delta c_{-3})e^{\lambda_{-3}a}\\
(\Delta c_{-2})e^{\lambda_{-2}a}\\
(\Delta c_{-1})e^{\lambda_{-1}a}\\
(\Delta c_{+1})e^{\lambda_{+1}b}\\
(\Delta c_{+2})e^{\lambda_{+2}b}
\end{bmatrix} 
=
\begin{bmatrix}
b_1 \\ b_2 \\ r \\ b_n \\ b_{n+1} \\ b_{n+2} \\ b_{n+3} \\ b_{n+4} \\ b_{n+5} \\ b_{n+6}
\end{bmatrix} \,,
\end{equation*}
where $M$ and $r$ denote the appropriate matrix and vector obtained from the discretisation of \eqref{eq:newtoniteration} at each point $\xi_2,\xi_3,\dots,\xi_{n-2}$ ($M$ is essentially the Jacobian, and $r$ is the negation of the most recent residual).
The remaining nine rows of this system are used to enforce the eight matching conditions \eqref{eq:asymptoticmatching} and the ``anchoring'' condition \eqref{eq:anchoringcondition}.
The specific ordering of these across the nine rows is unimportant, but it is convenient to enforce the matching condition for $\xi=a$ in the zeroth and first derivatives across rows 1 and 2, respectively, and the matching condition for $\xi=b$ in the zero'th and first derivatives across rows $n$ and $n+1$, respectively, i.e. as one might do for a more typical fourth order boundary value problem.
Notice that some scaling of the $\Delta c_j$ is factored into the column vector on the left hand side in order to improve the conditioning of the matrix (useful for very small $\delta$ or very large $\delta$ when the real part of the eigenvalues of the asymptotic solution can vary significantly in magnitude).

Accuracy of the solver is readily checked by: (a) varying the spatial resolution of the discretisation of the interval $[a,b]$; or (b) comparing the central part of solution when computed with larger intervals (and fixed step size).
With respect to both tests we observe extremely good agreement between solutions.
While we found the above approach to be sufficient with quadratic convergence rate, higher accuracy could be achieved by replacing the second order finite difference estimation of \eqref{eq:newtoniteration} with higher order finite difference stencils or any other appropriate collocation method.

\section{Comparison of various Boussinesq systems}\label{sec:comparison}

In this section, we employ the numerical method described in Section \ref{sec:computho} to compute heteroclinic orbits for the aforementioned systems. Our focus is on examining how the shape of these orbits varies with different values of $\delta > 0$ and $s > 1$. For the numerical computations, we used a convergence tolerance of $10^{-9}$ in Newton's method and discretized the domain $[a,b]$ using a uniform grid of $n+1$ nodes, where $n=2^4(b-a)$, assuming $a,b\in\mathbb{Z}$. 

In all the comparisons that follow, the same numerical value is used for the scaled dissipation parameters $\delta$ and $\bar{\delta}$. For the Bona--Smith system with $\gamma=1$ and the classical Boussinesq--Peregrine system this corresponds to the same viscosity $\nu$, while for the BBM--BBM system it corresponds to a viscosity smaller by a factor $\sqrt{2}$.

\begin{figure}[t]
\centering
\includegraphics[width=1\columnwidth]{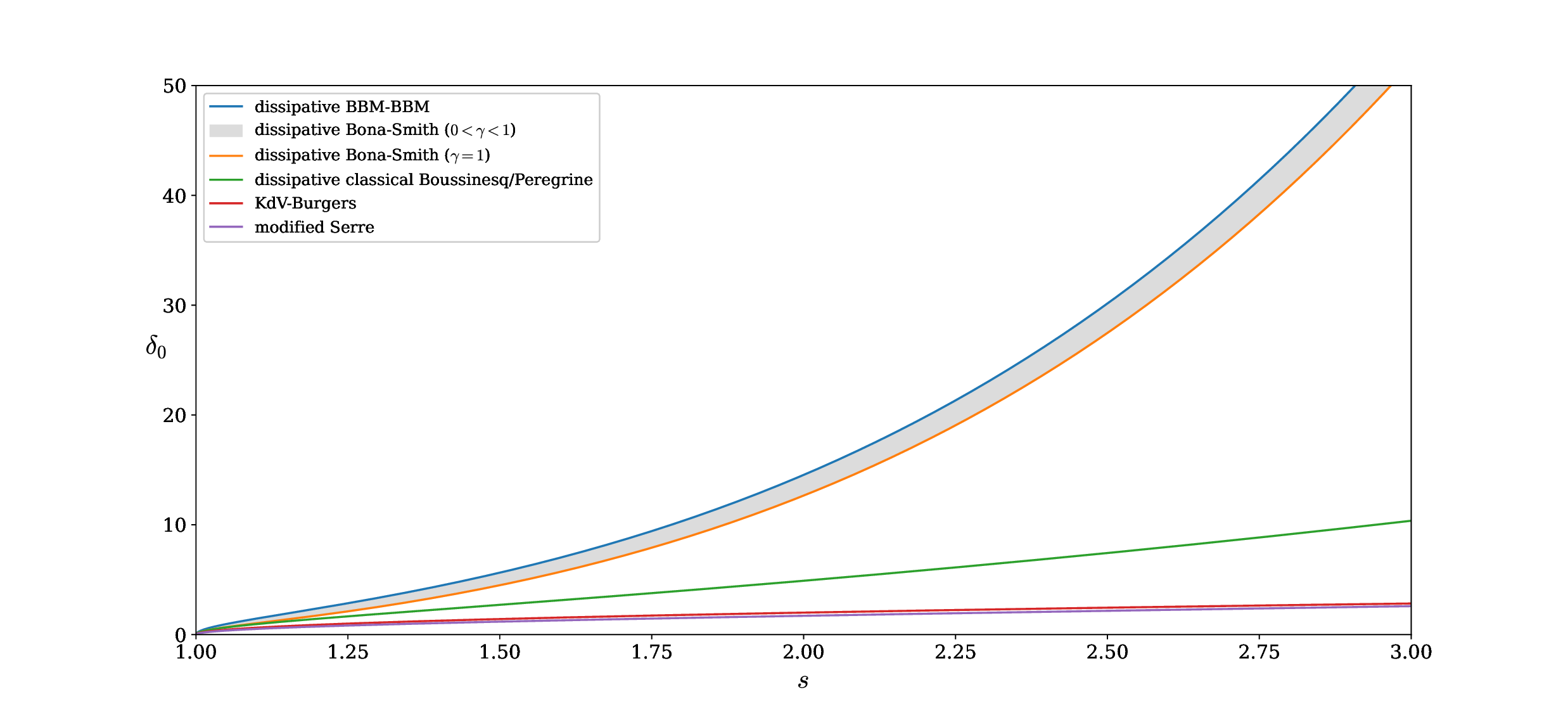}
\caption{Comparison of the bifurcation values $\delta_0$ for various dissipative-Boussinesq-type water wave models.}
\label{fig:comparedelta}
\end{figure}

First, we observe that the critical bifurcation values $\delta_0$ differ among the various Bona--Smith systems, depending on the value of $\gamma$. Furthermore, these critical values vary significantly across other Boussinesq-type models, such as the dissipative classical Boussinesq--Peregrine system \cite{BMT2022}, the modified Serre--Green--Naghdi system \cite{BMT2024}, and the KdV--Burgers equation \cite{BS1985}. Interestingly, the differences between the corresponding critical values $\delta_0$ across these models tend to diminish for small values of $s\approx 1$, but grow as $s$ increases away from $1$. This trend is illustrated in Figure \ref{fig:comparedelta}, where we plot $\delta_0$ and $\bar{\delta}_0$ as functions of $s$. The values of $\delta_0$ shown are the positive roots of the polynomial in (\ref{eq:condondelta}) and of $\tilde{S}(\delta)$, computed using the \texttt{numpy} function \texttt{roots}, for $s = 1.0001, 1.0002, \dots, 3$.

\begin{figure}[t]
\centering
\includegraphics[width=0.4\columnwidth]{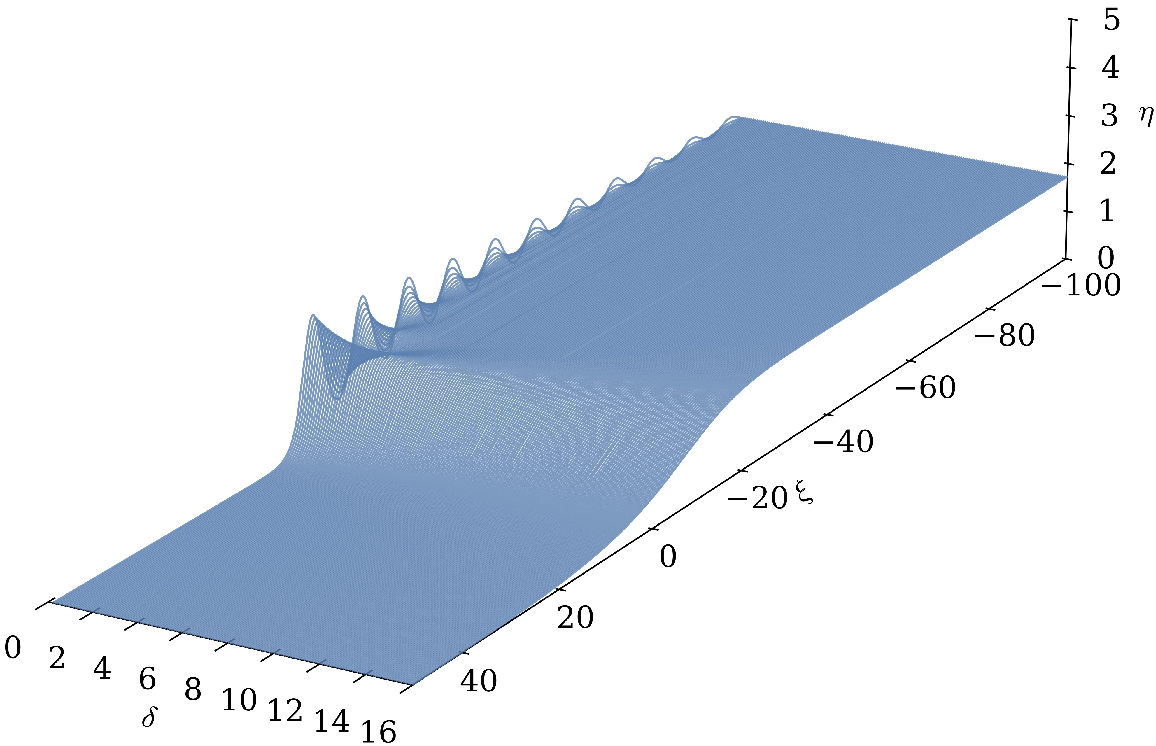}
\includegraphics[width=0.4\columnwidth]{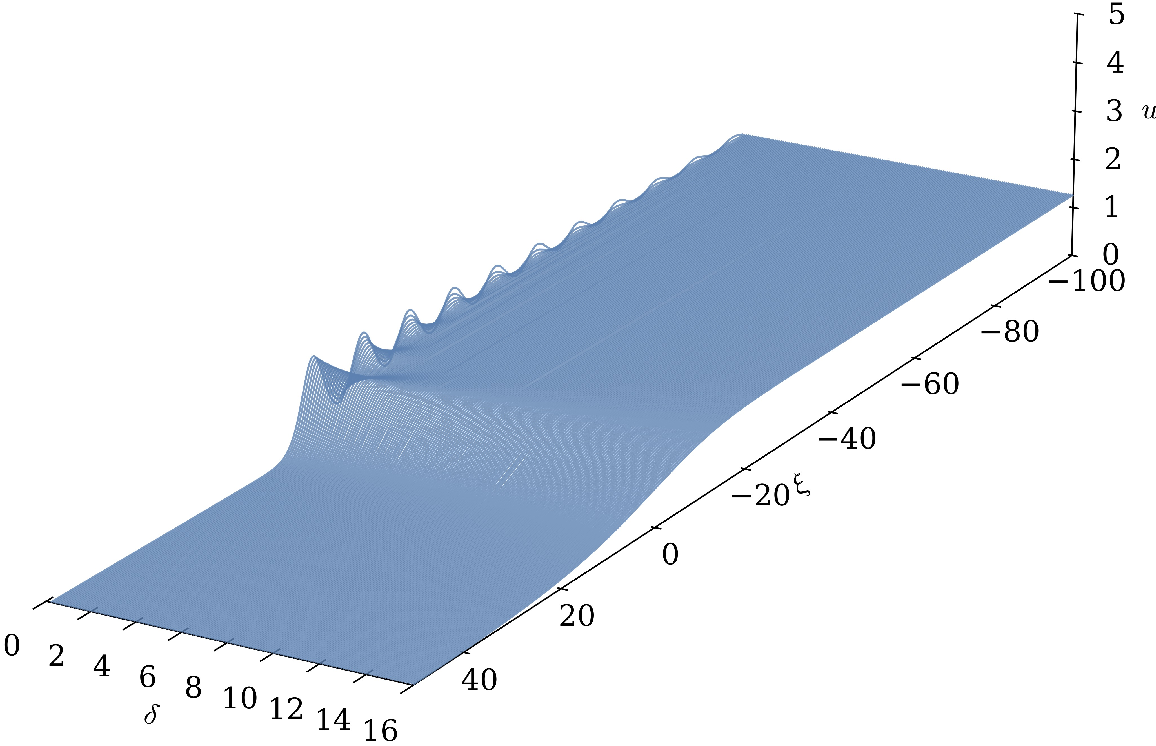}
\caption{Effect of dissipation parameter $\delta$ on traveling waves of the BBM--BBM system for $s=2$ and $\delta\in [0.4,16]$.}
\label{fig:comparedelta2}
\end{figure}
\begin{figure}[t]
\centering
\includegraphics[width=0.4\columnwidth]{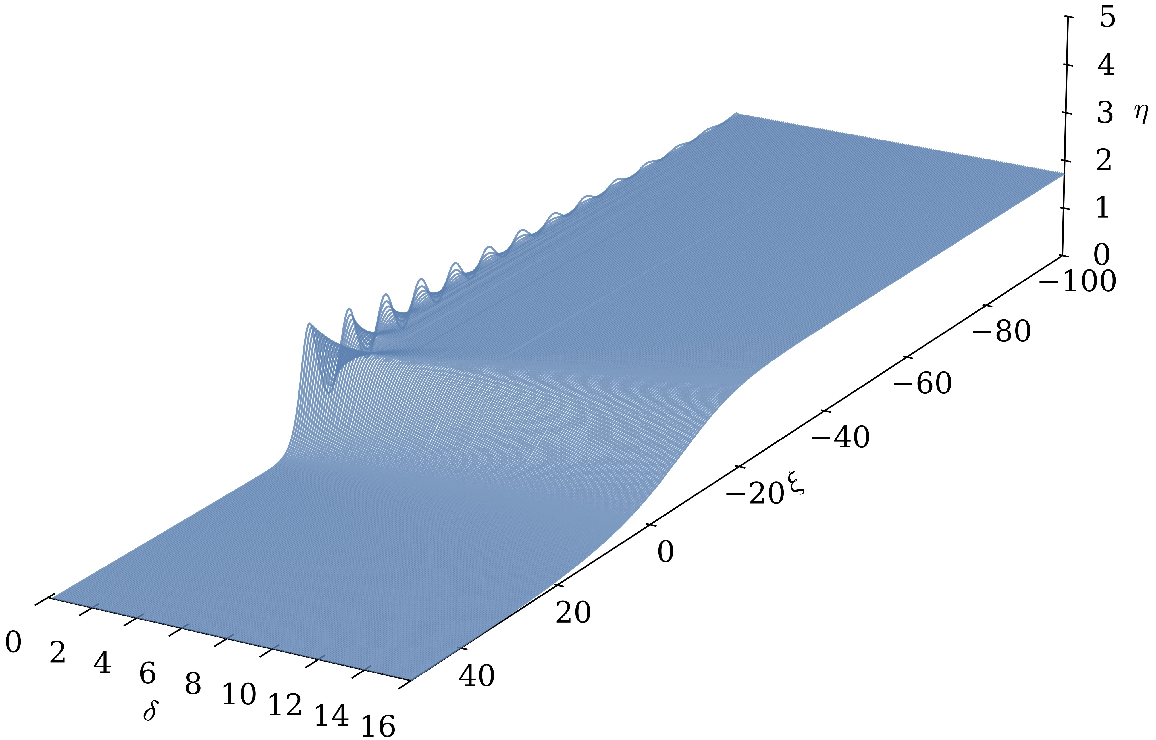}
\includegraphics[width=0.4\columnwidth]{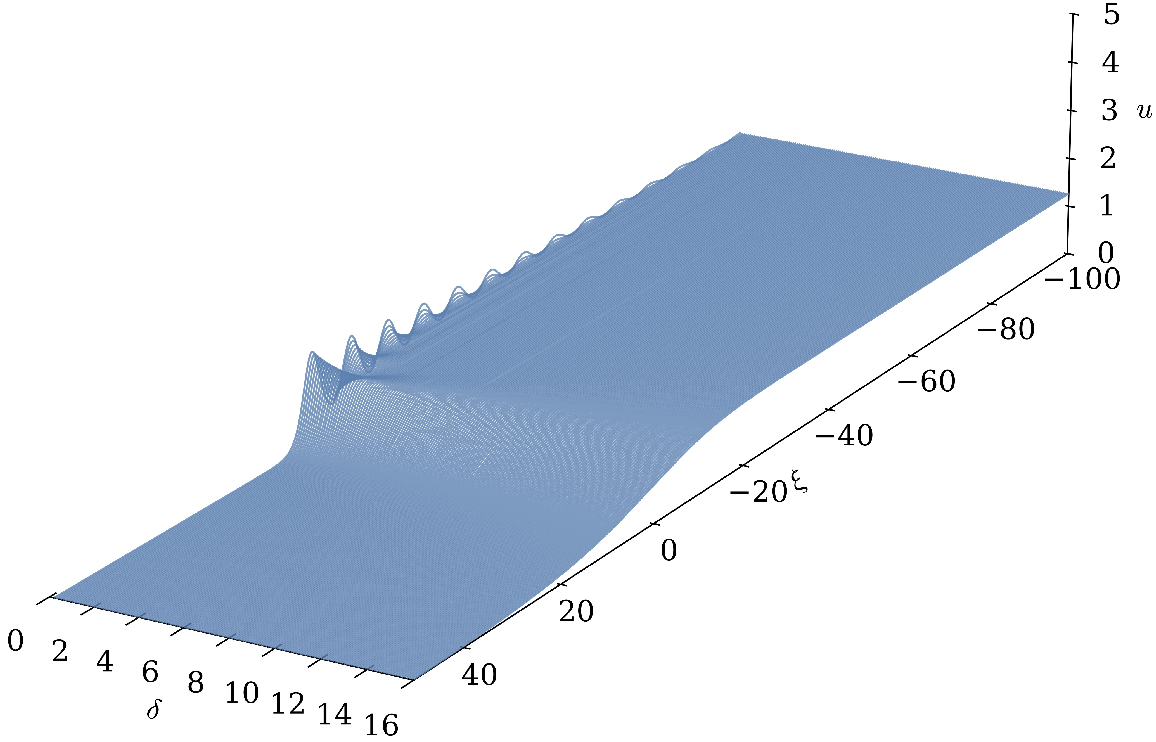}
\caption{Effect of dissipation parameter $\delta$ on traveling waves of the Bona--Smith system for $\gamma=1$, $s=2$ and $\delta\in [0.4,16]$.}
\label{fig:comparedelta2bs}
\end{figure}
\begin{figure}[t]
\centering
\includegraphics[width=0.4\columnwidth]{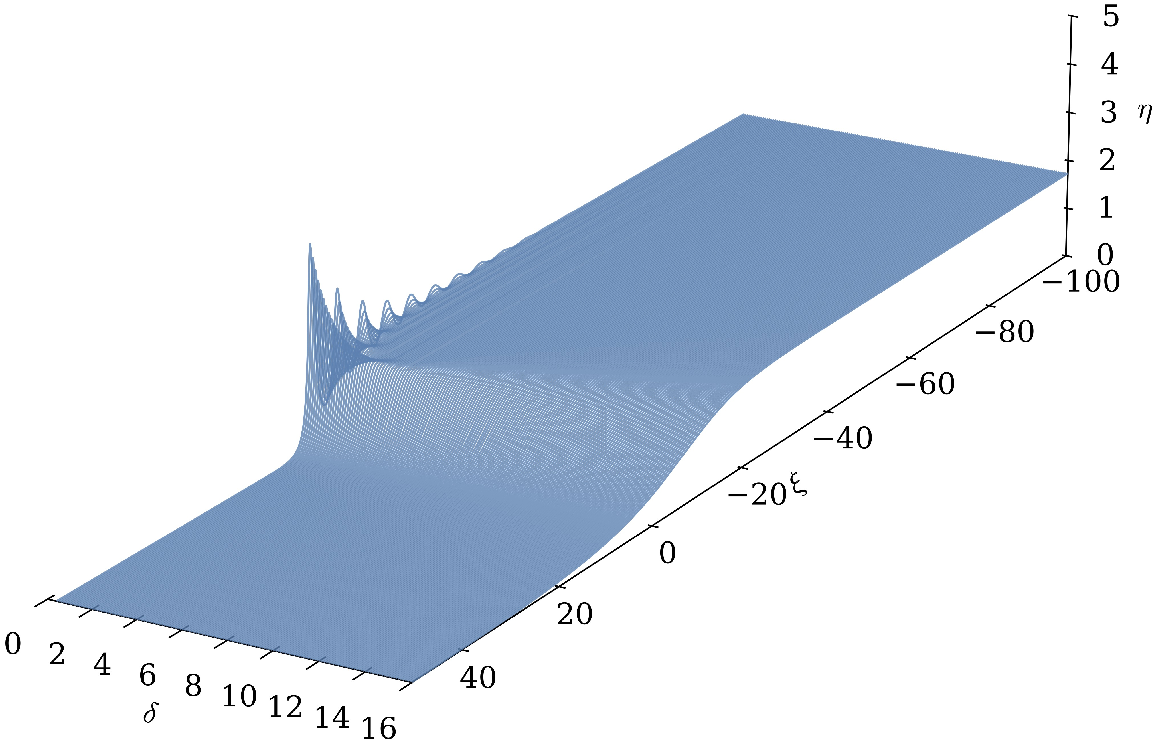}
\includegraphics[width=0.4\columnwidth]{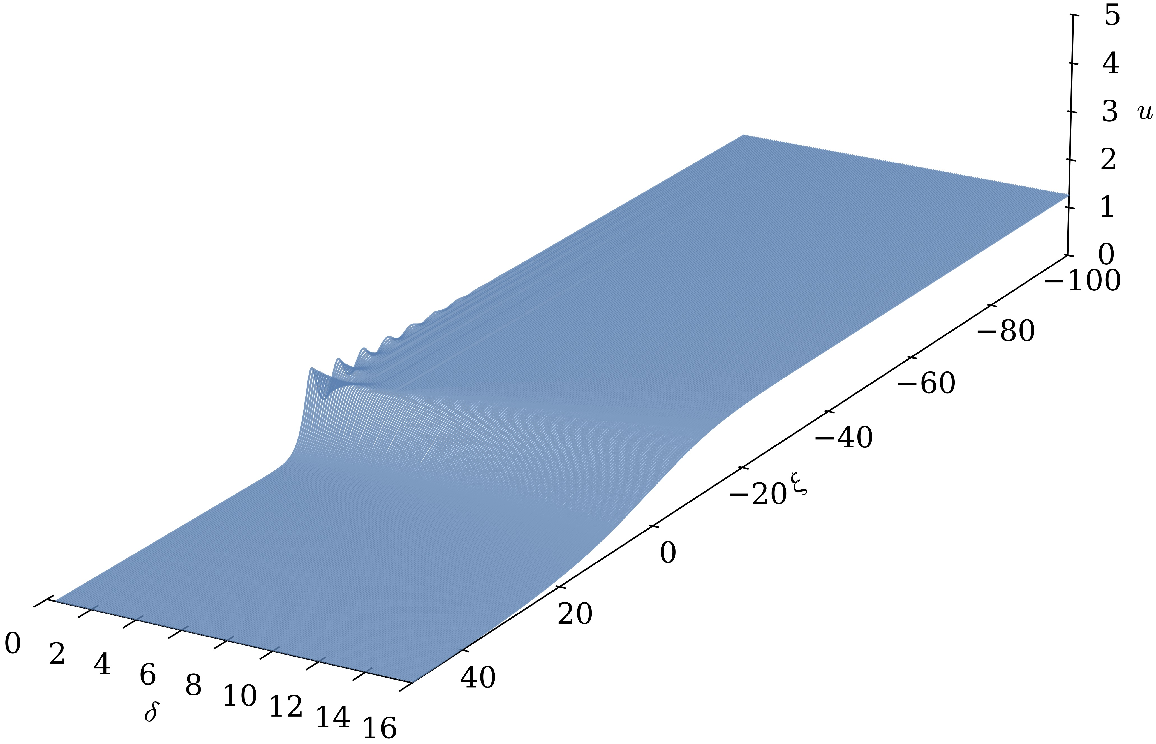}
\caption{Effect of dissipation parameter $\delta$ on traveling waves of the classical Boussinesq--Peregrine system for $s=2$ and $\delta\in [0.4,16]$}
\label{fig:comparedelta3}
\end{figure}

\begin{figure}[t]
\centering
\includegraphics[width=1\columnwidth]{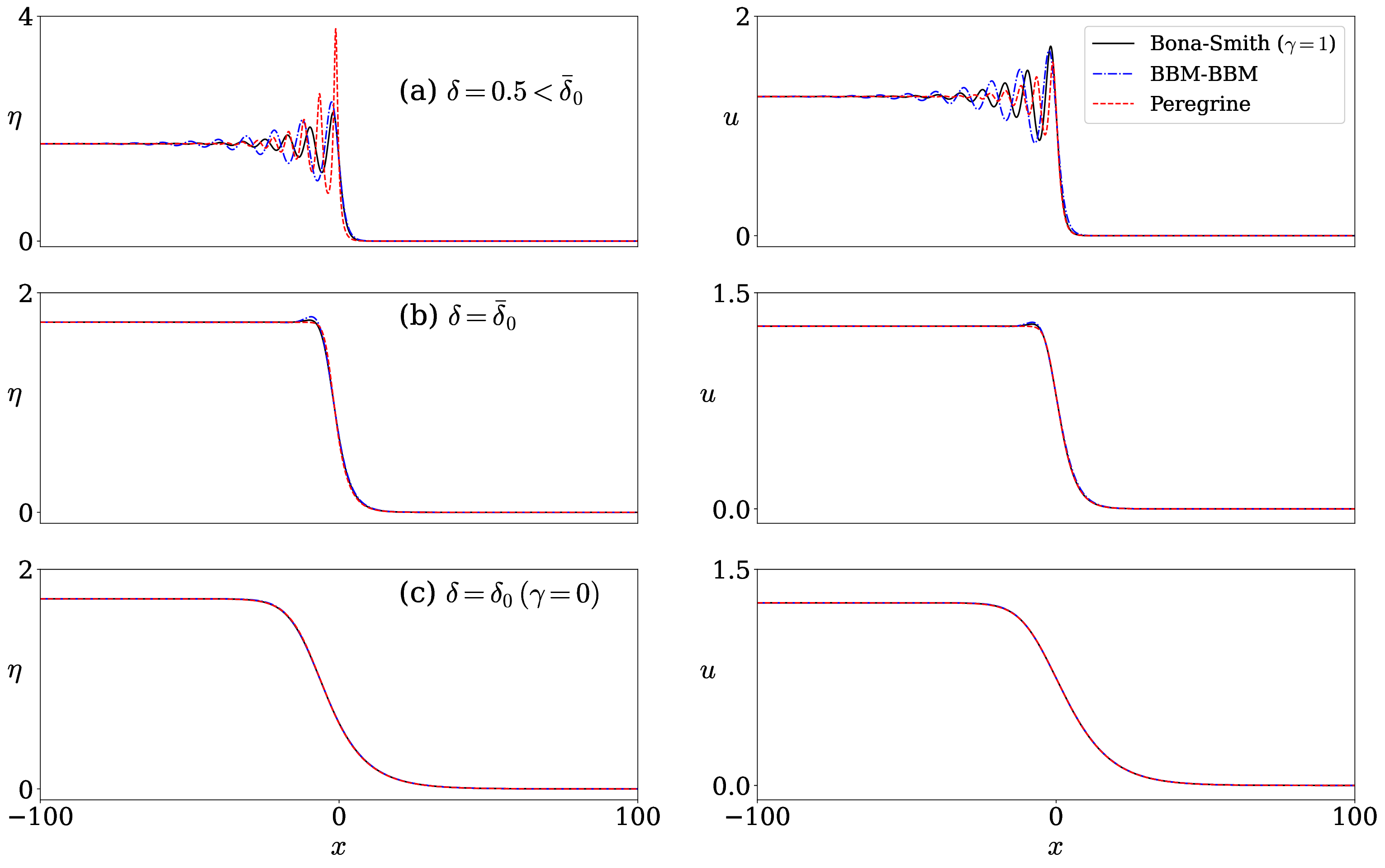}
\caption{Comparison of solutions for small and large values of $\delta$ and $s=2$.}
\label{fig:comparedelta4}
\end{figure}

We now investigate the shape of the solutions for the aforementioned systems across a range of $\delta$ values, fixing $s = 2$. Although similar qualitative behavior was observed for all values of $s > 1$ that we tested, we focus here on the representative case $s = 2$. This value leads to waves of moderate amplitude, with $\eta_0=\sqrt{3}$ in the scaled variables, for which the differences between the models are clearly visible. For values of $s$ close to $1$ the differences between the models become small, as indicated by the critical values $\delta_0$ shown in Figure \ref{fig:comparedelta}, while for larger values of $s$ the amplitude of the waves increases and the waves move away from the weakly nonlinear regime in which the Boussinesq systems are valid. Figures \ref{fig:comparedelta2}, \ref{fig:comparedelta2bs}, and \ref{fig:comparedelta3} display traveling wave solutions of systems (\ref{eq:sys1}), (\ref{eq:bsys1}), and (\ref{eq:cbscaled2}), respectively, for $\gamma = 1$ and $\delta \in [0.4, 16]$.

For reference, the critical values of $\delta$ for $s = 2$ are as follows: $\delta_0(\gamma=0) \approx 14.538370386372641$ for the BBM--BBM system, $\delta_0 \approx 12.65598071295166$ for the Bona--Smith system with $\gamma = 1$, and $\bar{\delta}_0 \approx 4.898979485566357$ for the classical Boussinesq system.

In these figures, we observe that for small values of $\delta$, the solutions of the Bona--Smith-type systems and the Peregrine system exhibit notable differences. As $\delta$ decreases to zero, the frequency and the magnitude of the oscillations increase in both models.  However, as $\delta$ increases, these decay rapidly, and their solutions become quite similar for values of $\delta$ only slightly greater than $0.5$. This convergence is consistent with the fact that these systems are asymptotically equivalent to each other. 

Nevertheless, even when $\delta = 0$, their solitary wave solutions differ in amplitude and spatial variation for the same value of $s$.
As before, the critical value $\delta_0$ plays a central role in determining the qualitative nature of the solution. For every value of $\gamma$ we considered, the solution becomes strictly decreasing when $\delta \geq \delta_0$, while regular oscillations appear for all $\delta < \delta_0$. Moreover, numerical computations confirmed that, in all cases, the maxima and minima of the solutions remain within the region $\mathcal{D}$.

\begin{figure}[t]
\centering
\includegraphics[width=1\columnwidth]{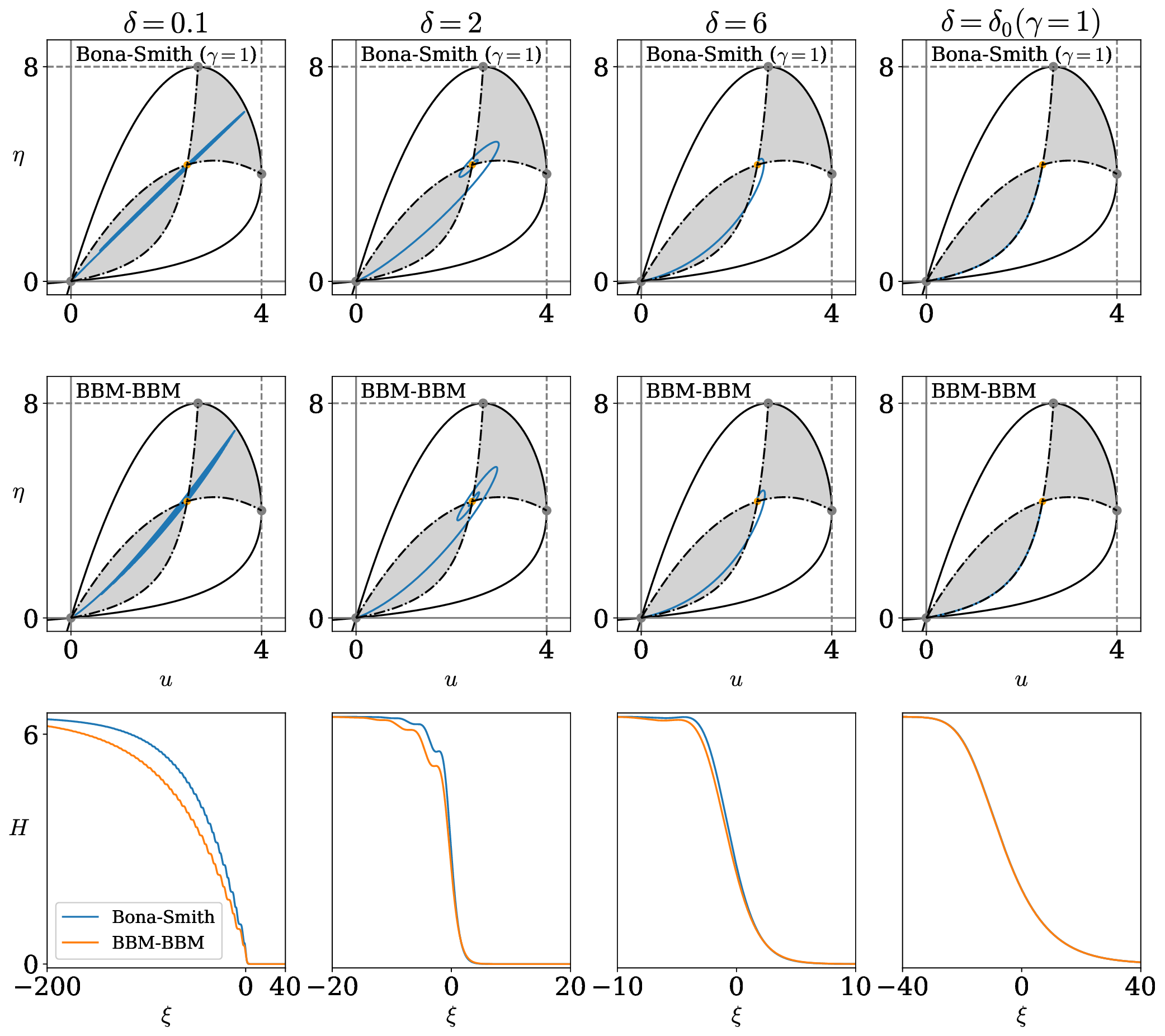}
\caption{Comparison of heteroclinic orbits of Bona--Smith ($\gamma=1$) and BBM--BBM systems in the physical space along with the corresponding graphs of $H(\bU)$ as a function of $\xi$ for $s=3$, $\delta=0.1$, $2$, $6$ and $\delta_0(\gamma=1)$}
\label{fig:comparephys}
\end{figure}

It is also worth emphasizing that, for a fixed value of $\delta$, the solutions of the various Boussinesq-type systems we examined are remarkably similar, particularly as $\delta$ becomes large. In fact, the solutions seem to converge to a common profile in the limit $\delta \to \infty$, corresponding to the regime where $f_u(\bu) = 0$. This behavior is illustrated in Figure \ref{fig:comparedelta4} and also in Figure \ref{fig:comparephys}, where the parameter $s = 3$ is used instead of $s = 2$.

We conclude that the parameter $\delta$ plays a key role in determining the shape of the system's solutions. For instance, in the case of the classical Boussinesq--Peregrine system, it is known that when $\bar{\delta} \geq \bar{\delta}_0$, the solution is strictly monotonic and is referred to as a regularized shock wave. In contrast, for $0<\bar{\delta} < \bar{\delta}_0$, the solution exhibits oscillatory behavior and is known as an oscillatory shock wave.

Figure \ref{fig:comparedelta4}(a) shows the solutions of the Bona--Smith system with $\gamma = 1$, the BBM--BBM system, and the classical Boussinesq--Peregrine system for $\delta = 0.5 < \bar{\delta}_0$, a value for which all solutions take the form of oscillatory shock waves. We observe notable differences not only in the amplitude of the oscillations but also in their frequency. As $\delta\to 0$, the leading pulse of the oscillatory shock wave is expected to approach a solitary wave of the corresponding conservative system with the same speed $s$. Since the solitary waves of these systems are different, substantial differences between the solutions are expected to persist for small but nonzero values of $\delta$. We note that heteroclinic orbits of this type do not exist for $\delta=0$, and thus no computations were performed in this case. These discrepancies, however, fall within the range allowed by their asymptotic equivalence, as all systems are formally consistent with the same underlying asymptotic regime.

In Figure \ref{fig:comparedelta4}(b), we present the solutions of the same systems for $\delta = \bar{\delta}_0$. At this value, only the solution of the classical Boussinesq--Peregrine system takes the form of a regularized shock wave, that is, a strictly decreasing profile, while the solutions of the other systems remain in the oscillatory shock wave regime. Despite the significant differences in the critical bifurcation values $\delta_0$ across these models, the resulting wave profiles are remarkably similar. This close agreement highlights the physical relevance and consistency of all these models in capturing the same underlying wave dynamics. It is noted that for the computed solutions of the Bona--Smith system, $H(\bU)-f(\bu_0)<-10^{-3}<0$ for all values $\xi$ we computed.

Finally, in Figure \ref{fig:comparedelta4}(c), we present the solutions of these systems for $\delta = \delta_0(\gamma = 0)$, a value for which all models are expected to support regularized shock waves. The resulting wave profiles appear virtually identical, though minor discrepancies become apparent upon closer examination or through detailed numerical comparison.

Figure \ref{fig:comparephys} shows the orbits obtained for $\delta = 0.1$, $2$, $6$, and for $\delta_0(\gamma = 1)$, with $s = 3$. The first row displays the orbits corresponding to the Bona--Smith system with $\gamma = 1$, while the second row presents the corresponding orbits for the BBM--BBM system. As $\delta$ decreases, the orbits tend to approach the arc $P_1P_2$. In our computations, we did not observe any orbit intersecting this arc; however, it is expected that such an intersection would occur in the limiting case $\delta = 0$, where the orbit becomes homoclinic and attains its maximum value on the arc $P_1P_2$ \cite{Toland1986}. In Figure \ref{fig:comparephys}, we do not include solutions from the classical Boussinesq--Peregrine system because all its orbits lie on the curve defined by $f_u(\bu) = 0$. 

Although the points $(u,\eta)$ that correspond to extrema of the solutions of the Bona--Smith system do not necessarily occur in the same region as those of the BBM--BBM system, in Figure \ref{fig:comparephys} we present the region where extrema of the solutions of the BBM--BBM system are most likely to appear. It seems that both models behave in a very similar manner, at least for sufficiently large values of $\xi$ with extrema in both cases to occur in the same region. For very negative values of $\xi$ it is hard to verify the location of the extrema due to floating point arithmetic.

These plots also suggest that the orbits in $\mathbb{R}^4$ lie on a three-dimensional manifold that becomes nearly orthogonal to the physical plane as $\delta \to 0$, resulting in a symmetric homoclinic orbit. Conversely, as $\delta$ increases, this manifold appears to flatten above the physical space. Also, for larger values of $\delta$, the solutions are confined to regions progressively further from the arc $P_1P_2$, and in all cases, the orbits appear to converge toward a solution lying on the line $f_u(\bu) = 0$ as $\delta$ increases. To verify this, consider the change the variable $\xi\to \xi/\delta$. Then, the first equation of (\ref{eq:sys2}) becomes
$$-\frac{s}{\delta^2}\eta''+f_u(\bu)=0\ ,$$
which implies that the solution should satisfy $f_u(\bu)\sim 0$ as $\delta\to\infty$ for all $s>1$ and $\gamma\in [0,1]$. As it is suggested by the computed orbits, the solution converges to such a regime very fast.

The last row of Figure \ref{fig:comparephys} compares the graphs of $H(\bU)$ for the Bona--Smith system ($\gamma = 1$) and the BBM--BBM system ($\gamma = 0$). We observe that whenever $u'$ or $\eta'$ vanish, the derivative $\frac{d}{d\xi} H$ also becomes zero, in accordance with equation (\ref{eq:potender}). Moreover, in all cases, the values of $H(\bU)$ satisfy $0 < H(\bU) < f(\bu_0)$. Note that the $\xi$-axis scales differ between these graphs to better illustrate the details of $H$. As $\delta$ decreases, the value of $H$ at the inflection point of $u$ at $\xi = 0$ approaches the origin, consistent with the expectation that $H$ vanishes in the limit $\delta = 0$, where the systems become Hamiltonian with $H(\bU) = 0$. These findings align with the theoretical discussion of Section \ref{sec:existence}.

\begin{figure}[t]
\centering
\includegraphics[width=1\columnwidth]{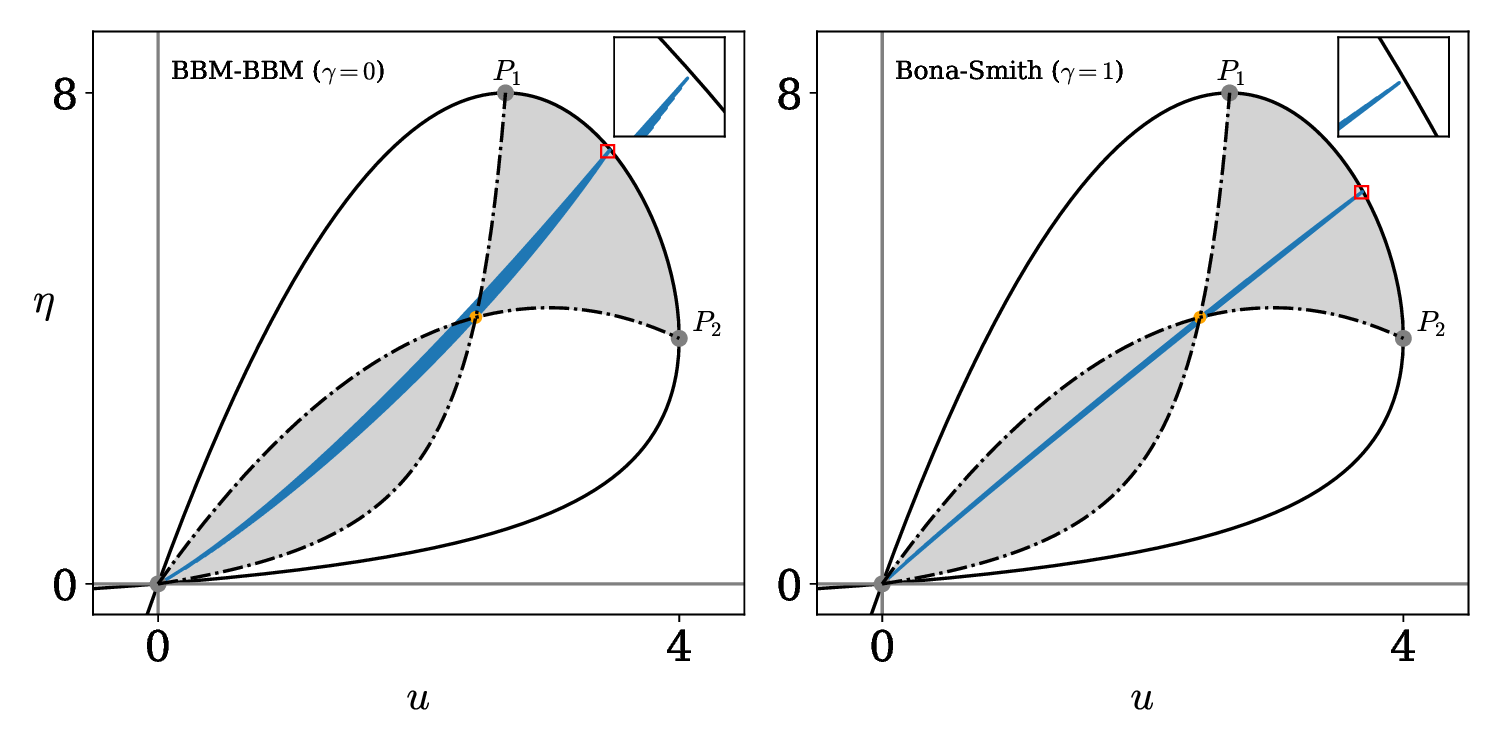}
\caption{Oscillatory shock waves for $s=3$ and $\delta=0.01$}
\label{fig:comparephys2}
\end{figure}

As briefly discussed in Section \ref{sec:computho}, the convergence of our numerical method becomes increasingly challenging as $\delta$ approaches very small values, particularly for larger values of $s$. Our approach is to first calculate an initial orbit for a large value of $\delta$ and  $s$ close to $1$, and then perform continuation by gradually decreasing $\delta$ and increasing $s$ using small steps (usually $O(10^{-3})$) toward desired values. At each step, the previously computed orbit serves as the initial guess for Newton's method, ensuring convergence.  

The smallest value of $\delta$ we successfully tested for $s=3$, $\gamma=0$ and $\gamma=1$ was $\delta = 0.01$. To compute the heteroclinic orbits for such values of $\delta$ and $s$, we first obtained the orbits at $s = 1.01$ and $\delta = 0.1$. We then applied continuation in $\delta$ with a step size of $d\delta = -0.01$ to reach $\delta = 0.01$. Following this, continuation with respect to $s$ using a step size of $ds = 0.0001$ allowed us to compute the orbits for the desired values of $s$ and $\delta$. Figure \ref{fig:comparephys2} presents the computed heteroclinic orbits in physical space and confirms that the orbit does not intersect the curve $\mathcal{C}$. The figure also includes a magnified view of the orbit near the curve $\mathcal{C}$, highlighting the inclusion of the orbits in $\mathcal{D}$ in greater detail. However, we cannot exclude the possibility that, for smaller values of $\delta$, the orbits of this system may exit the region $\mathcal{D}$. In such a case, we would expect that only a very limited portion of the orbit would leave the open arc $P_1P_2$. It is worth noting that, for a heteroclinic orbit of interest to exit from $\mathcal{C} \setminus \overline{P_1P_2}$, the functional $H(\bU)$ must become negative inside $\mathcal{D}$; otherwise, it cannot re-enter $\mathcal{D}$ and converge to the origin.

We note that for the traveling waves computed (but not shown) for verification purposes in this study, we employed sufficiently large intervals, (varying according to the specific case), so that the solution values at the boundaries effectively coincide with the corresponding equilibria. Using these computed solutions as initial conditions for system (\ref{eq:sysbs2}), we further verified numerically, via the modified finite element method of \cite{MM2023}, that the waves propagate without change in shape (within numerical accuracy). 

\begin{figure}[t]
\centering
\includegraphics[width=1\columnwidth]{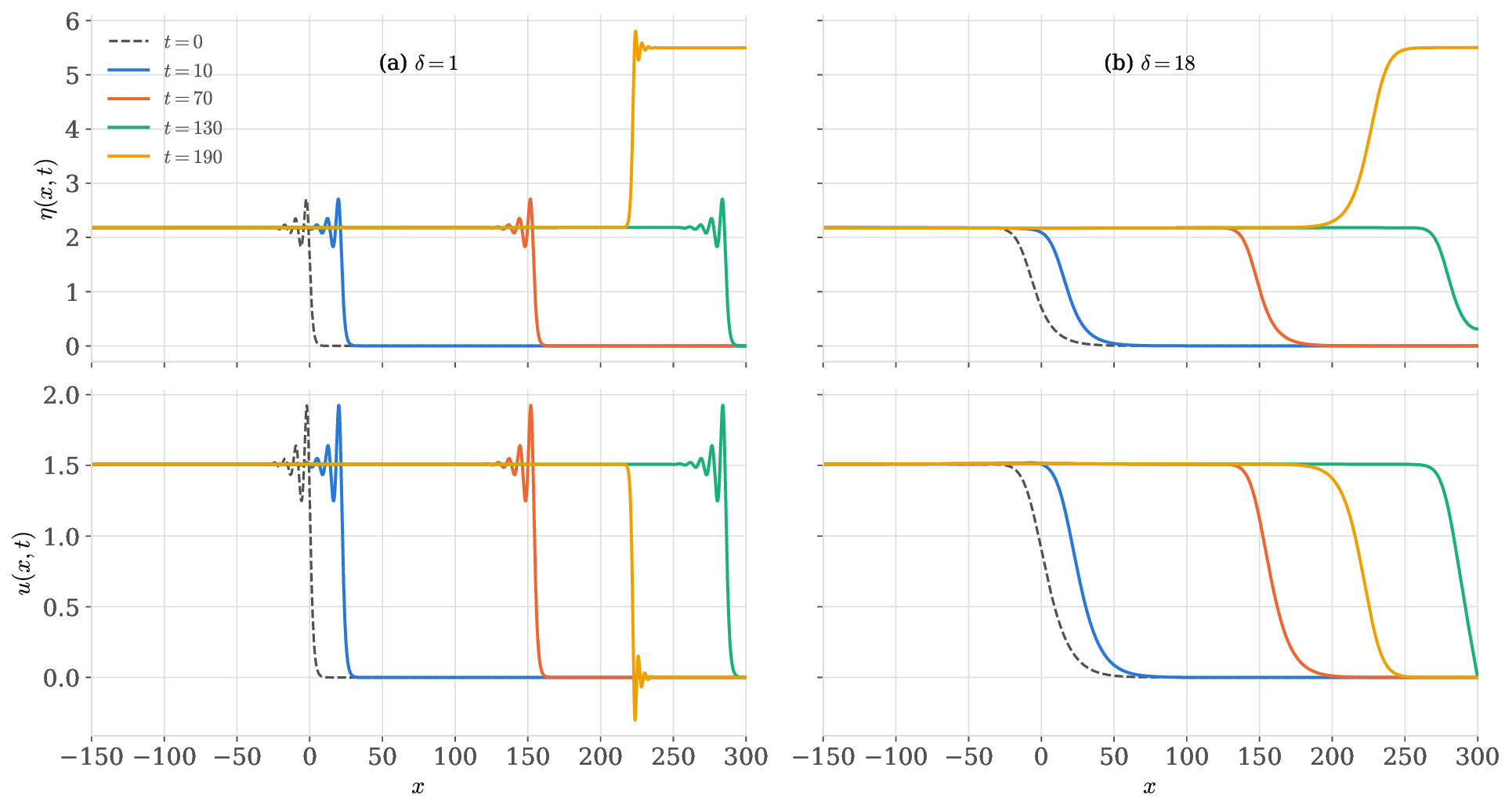}
\caption{Propagation and reflection of traveling waves with $s=2.2$ for the Bona--Smith system ($\gamma=1$) and dissipation parameters (a) $\delta=1$ and (b) $\delta=18$}
\label{fig:femstab}
\end{figure}

For instance, Figure \ref{fig:femstab} shows the propagation of an oscillatory shock wave with $s=2.2$ for the Bona--Smith system with $\gamma=1$ and $\delta=1$ and of a regularized shock wave with the same $\gamma=1$ but with $\delta=18$ in the interval $(a,b)=(-150,300)$. The value $\delta=18$ is larger than the critical value $\delta_0(\gamma=1)\approx 17.6253$, and thus the equilibrium $\bU_0$ is a saddle and the traveling wave is strictly monotone. As an initial conditions, we employed the corresponding orthogonal projections onto the finite element space of a numerically generated traveling waves obtained by the previously described method. In this experiment, the finite element space was chosen to be the standard Lagrange space of piecewise linear elements with spatial mesh length $\Delta x=0.1$. Time integration was carried out using the classical explicit fourth-order Runge--Kutta scheme with four stages and with a time step $\Delta t=0.1$. The results displayed in the figure suggest that these traveling waves are stable. 

Moreover, because we imposed reflection boundary conditions at the right endpoint of the interval $(a,b)$, namely $\eta_x(a)=\eta_x(b)=0$, $u(a)=u_0$, $u(b)=0$, \cite{ADM2009}, the traveling waves were reflected at the wall and then propagated to the left. The reflection of waves in such systems is equivalent to the head-on collision of two identical, counter-propagating waves \cite{ADM2010}, and shows that these waves remain stable even under large perturbations. Comparable conclusions can be drawn for analogous experiments involving other Bona--Smith systems, including the BBM--BBM system. Further theoretical and numerical analysis of the solutions to such models will be the subject of future work.

\begin{figure}[t]
\centering
\includegraphics[width=1\columnwidth]{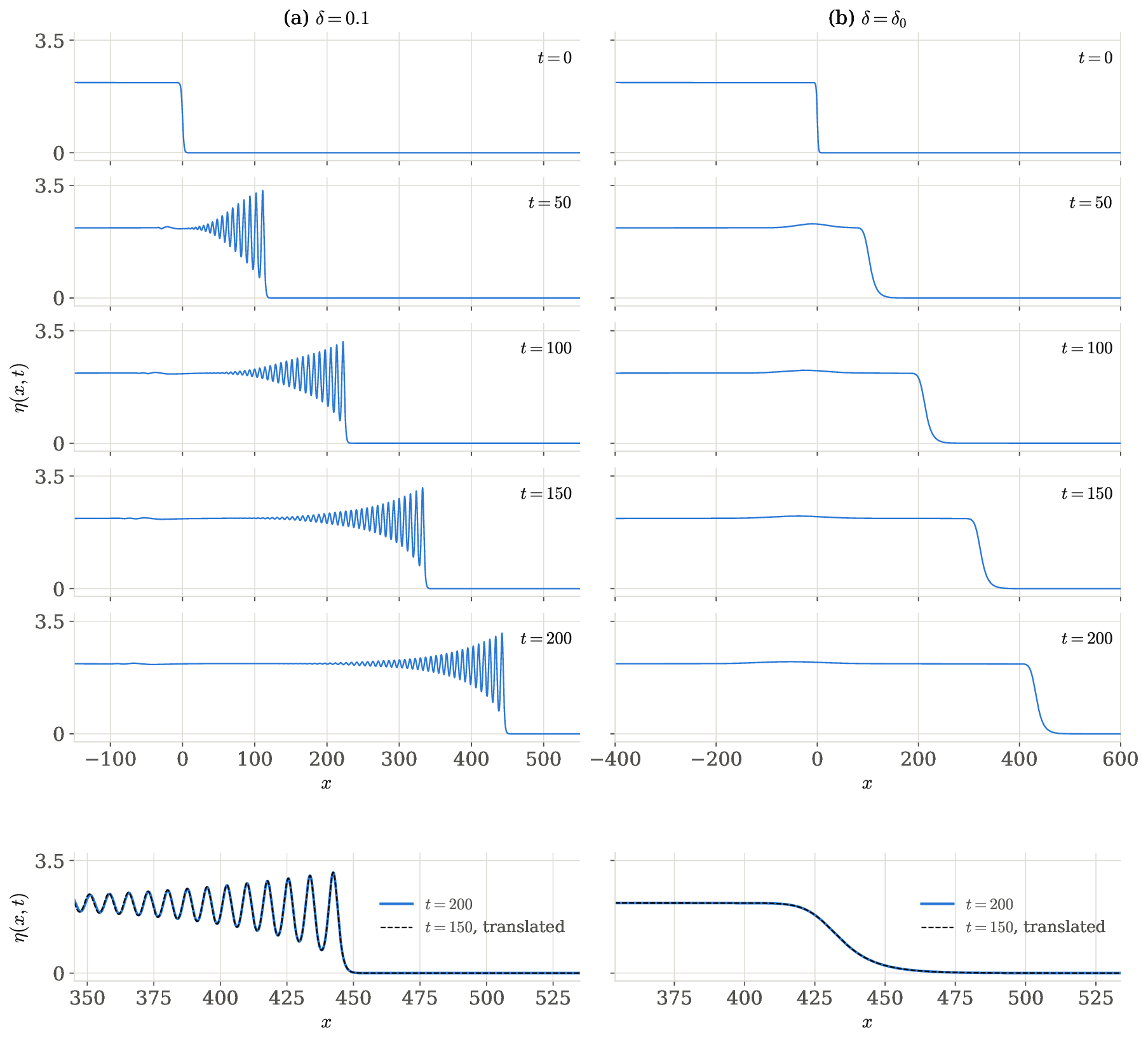}
\caption{Generation of traveling waves of the Bona--Smith system ($\gamma=1$) from the smoothed step (\ref{eq:step}) with $s=2.2$: (a) oscillatory shock wave for $\delta=0.1$ and (b) regularized shock wave for $\delta=\delta_0$. The bottom panels show the solution at $t=200$ around the front and the solution at $t=150$ translated so that the two fronts coincide.}
\label{fig:generation}
\end{figure}

Finally, we examine whether these traveling waves are generated by initial conditions that are not traveling waves. We consider the Bona--Smith system with $\gamma=1$ and the smoothed step
\begin{equation}\label{eq:step}
\eta(x,0)=\frac{\eta_0}{2}\left(1-\tanh\frac{x}{w}\right)\ ,\qquad u(x,0)=\frac{u_0}{2}\left(1-\tanh\frac{x}{w}\right)\ ,
\end{equation}
with $w=2$, where $\bu_0=(u_0,\eta_0)$ is the equilibrium that corresponds to the speed $s=2.2$. We solved the system up to $T=200$ with the same numerical method and the same values of $\Delta x$ and $\Delta t$, for $\delta=0.1$ in the interval $(-150,550)$ and for the critical value $\delta=\delta_0(\gamma=1)\approx 17.6253$ in the interval $(-400,600)$. The results are presented in Figure \ref{fig:generation}. In both cases the step evolves into a front that travels to the right. The speed of the front, measured at the point where $\eta=\eta_0/2$, tends to $s$, and between $t=150$ and $t=200$ it is equal to $2.2000$. For $\delta=0.1$ an oscillatory shock wave is generated, while for $\delta=\delta_0$ the front is monotone and a regularized shock wave is generated. The bottom panels of Figure \ref{fig:generation} show the solution at $t=200$ around the front, together with the solution at $t=150$ translated so that the two fronts coincide. For $\delta=\delta_0$ the two profiles differ by less than $4\times 10^{-4}$ in the region shown. For $\delta=0.1$ they agree near the front and differ slightly further behind, where the oscillatory tail is still developing. Since the initial condition is not a traveling wave, a small disturbance of the left state also remains far behind the front and decays slowly.

\section{Conclusions}\label{sec:conclusions}

In this work, we investigated the influence of dissipation on traveling wave solutions of Bona--Smith-type systems arising in water wave theory. The analysis is particularly intricate due to the four-dimensional nature of the underlying dynamical systems. It is known that in the absence of dissipation, such systems admit symmetric homoclinic orbits as well as other multi-pulse nonphysical homoclinic orbits. As dissipation is introduced, the symmetric homoclinic orbits appear to deform into heteroclinic connections between two distinct equilibria. One equilibrium remains a saddle throughout, while the other (from the unperturbed saddle-center) undergoes a sequence of bifurcations, becoming first a saddle-focus and eventually a saddle as the dissipation parameter increases. Despite these changes, the resulting heteroclinic orbits seem to remain confined within a leaf-shaped domain of the physical space, whose geometry is independent of the dissipation parameter.

We showed that the qualitative structure of the heteroclinic orbits is governed by the magnitude of dissipation. For dissipation values below a critical threshold, the orbits exhibit an oscillatory shock profile, characterized by a sharp leading edge followed by decaying oscillations. In contrast, when the dissipation exceeds or meets the critical value, the solution becomes monotonic, resembling a regularized shock wave with no oscillatory tail. Contrary to other Boussinesq models, the local extrema of the oscillatory shock waves do not occur simultaneously resulting to a non-monotone Hamiltonian. 

To compute these heteroclinic orbits numerically, we developed an efficient second-order accurate finite difference method that circumvents the issue of translational invariance by anchoring an inflection point on the $y$-axis. The method also incorporates asymptotic estimates to enhance convergence.

We conclude by comparing the heteroclinic orbits arising in several Boussinesq-type systems, including the classical Boussinesq/Peregrine system. Our results suggest that such heteroclinic orbits are structurally stable. However, one can observe substantial differences in the wave profiles generated by each model but for large values of the dissipation parameter, the solutions tend to converge toward a single limiting orbit.

\appendix
\section{Proof of Proposition \ref{prop:eigenvalues}}\label{app:linearization}

In this appendix we present the detailed analysis of the roots of the characteristic polynomial $P_2(\lambda)$ given in (\ref{eq:P2poly}), which leads to Proposition \ref{prop:eigenvalues}. We use the notation of Section \ref{sec:dynamics}, and in particular the quantities $A$ and $q$ defined there.

Using Descartes' rule of signs together with the facts that $P_2(0) < 0$ and $\lim_{\lambda \to \pm \infty} P_2(\lambda) = +\infty$, we conclude that the polynomial $P_2$ always has exactly one negative real root and either one or three positive real roots, counted with multiplicity. Thus, besides the negative root $\lambda_r^-$ and the positive root $\lambda_r^+$, the remaining two roots are either a pair of complex conjugates or two positive real roots, depending on the value of $\delta > 0$.

Assuming that the eigenvalues are of the form $\lambda=\lambda_0+\delta \lambda_1+O(\delta^2)$, we find that the eigenvalues of $P_2(\lambda)$ as $\delta\to 0$ are
$$\lambda_r^\pm\sim \frac{\delta}{4s}\pm \rho_r+O(\delta^2), \quad \lambda_i^\pm\sim \frac{\delta}{4s}\pm \I \rho_i+O(\delta^2)\ , $$
indicating that the real part of the complex eigenvalues is increasing with $\delta$ but is less than the positive eigenvalue for at least small values of $\delta$. 

Note that $P_2(\lambda)$ has two inflection points for all $\delta>0$ and $s>1$. These are located at
\begin{equation}\label{eq:inflpts}
p_\pm= \frac{1}{24}\left(\frac{6\delta}{s}\pm\sqrt{36\frac{\delta^2}{s^2}+192\frac{s-u_0}{s}}\right)=\frac{\delta\pm\sqrt{\delta^2-\frac{8}{3}s\left(s-\sqrt{s^2+8}\right)}}{4s}\ .
\end{equation}
So $p_-<0$ and $p_+ >0$. Moreover, using the equation $P_2''(p_+)=0$ to eliminate the powers of $p_+$ that are higher than one, we obtain
$$P_2(p_+)=-\frac{\delta^3}{8s^3}\,p_+-\frac{A\delta^2}{12s^2}-\frac{5A^2}{9}-q<0\ .$$
Since $P_2(\lambda)$ is convex in $(p_+,+\infty)$ and $P_2(+\infty)=+\infty$, it has exactly one root in this interval and this root is simple. For $\delta=0$ this root is $\rho_r$. Therefore, it is the root $\lambda_r^+$ for all $\delta>0$, and we have $\lambda_r^+>p_+$, while all the other real roots of $P_2(\lambda)$ are less than $p_+$.

Applying Descartes' rule of signs to the derivative 
$$P_2'(\lambda)=4\lambda^3-3\frac{\delta}{s}\lambda^2-4\frac{s-u_0}{s}\lambda+\frac{\delta}{s}\frac{s-u_0}{s}\ ,$$
and also to $P_2'(-\lambda)$, and using the fact that for all $\delta>0$ and $s>1$, we have $P_2'(+\infty)=+\infty$ and 
$$P_2'(p_+)=-\frac{9 \delta^3 + (3 \delta^2  + 16 s(s - u_0))
  \sqrt{9 \delta^2 + 48 s (s - u_0)} }{72 s^3}<0\ ,$$ 
we deduce that $P_2'(\lambda)$ has one negative and two positive real roots for all $\delta>0$ and $s>1$. These will be located between $(\lambda_r^-,\lambda_r^+)$ due to Gauss--Lucas theorem. In what follows, we will denote the roots of $P_2'(\lambda)$ by $x_1<0$, $x_2>0$, $x_3>0$. 

The existence of a unique double eigenvalue indicates that the complex eigenvalues of $P_2(\lambda)$ can collapse onto a double real eigenvalue. This can happen only if the greatest common divisor between $P_2(\lambda)$ and $P_2'(\lambda)$ is a linear polynomial. However, using Euclid's algorithm, the last remainder of Euclid's algorithm is the constant
\begin{equation}\label{eq:gcd} \text{\footnotesize $\frac{
(3 \delta^2 + 16 A s^2)^2 (4 A^3 \delta^6 + 
     \delta^4 (A^2 - 3 q) (13 A^2 + 9 q) s^2 + 
     32 A \delta^2 (A^2 - 3 q) (A^2 + q) s^4 - 
     256 q (A^2 + q)^2 s^6)}{64 (3 A \delta^4 s + 
     2 \delta^2 (7 A^2 + 3 q) s^3 + 32 A (A^2 + q) s^5)^2}$}\ ,
\end{equation}
where 
\begin{equation}\label{eq:Aformula1}
    A=\frac{s-u_0}{s}>0\quad\text{and}\quad q=\frac{1+\eta_0-(s-u_0)^2}{s^2}>0\ .
\end{equation}
The only case where the greatest common divisor between $P_2(\lambda)$ and $P_2'(\lambda)$ is a linear polynomial is when the numerator in (\ref{eq:gcd}) is zero. Specifically, for given $s>1$, the numerator becomes zero when
\begin{equation}\label{eq:condondelta}
4 A^3 \delta^6 + 
     (A^2 - 3 q) (13 A^2 + 9 q) s^2 \delta^4 + 
     32 A  (A^2 - 3 q) (A^2 + q) s^4 \delta^2 - 
     256 q (A^2 + q)^2 s^6=0\ ,
\end{equation}
or equivalently by setting $x=\frac{\delta^2}{s^2}$ we have
\begin{equation}\label{eq:condondelta2}
4 A^3 x^3 + 
      (A^2 - 3 q) (13 A^2 + 9 q) x^2 + 
     32 A (A^2 - 3 q) (A^2 + q) x - 
     256 q (A^2 + q)^2=0\ .
\end{equation}
Since $A^2-3q$ is either positive or negative, applying Descartes' rule of signs in both cases, we deduce that there is a unique positive solution of the equation (\ref{eq:condondelta}) for $\delta^2$. We denote the positive root of (\ref{eq:condondelta}) by $\delta_0$. 

The root $\delta_0^2$ can be computed approximately by numerical means or exactly using Mathematica\textsuperscript{\tiny\textregistered}, which is
\begin{equation}\label{eq:critdelta0}
\delta_0^2=\frac{a_1}{12 A^3}\ ,
\end{equation}
where
$$\begin{aligned}
a_1 &=(27 q^2+30 A^2 q -13 A^4 ) s^2+ a_4+ \frac{s^4
   \left(-215 A^8-12 A^6 q+1350 A^4 q^2+1620 A^2 q^3+729 q^4\right)}{a_4}
   \\
a_2 &=(5291 A^{12} +38250 A^{10} q +85725 A^8 q^2 +102060 A^6 q^3+91125 A^4 q^4 +65610 A^2 q^5) s^6\\
a_3 &= \left(7 A^6+25 A^4 q+45 A^2 q^2+27 q^3\right)^3 \\
a_4 &= \sqrt[3]{a_2+192 \sqrt{3} \sqrt{A^6 s^{12} a_3}+19683 q^6 s^6}\ .
\end{aligned}
$$
Since $P_2'(\lambda)$ has always one negative and two positive real roots, we have that $P_2(\lambda)$ has two complex roots for $0<\delta<\delta_0$ and no complex roots for $\delta\geq \delta_0$.
 
\begin{figure}[t]
\centering
\includegraphics[width=0.8\columnwidth]{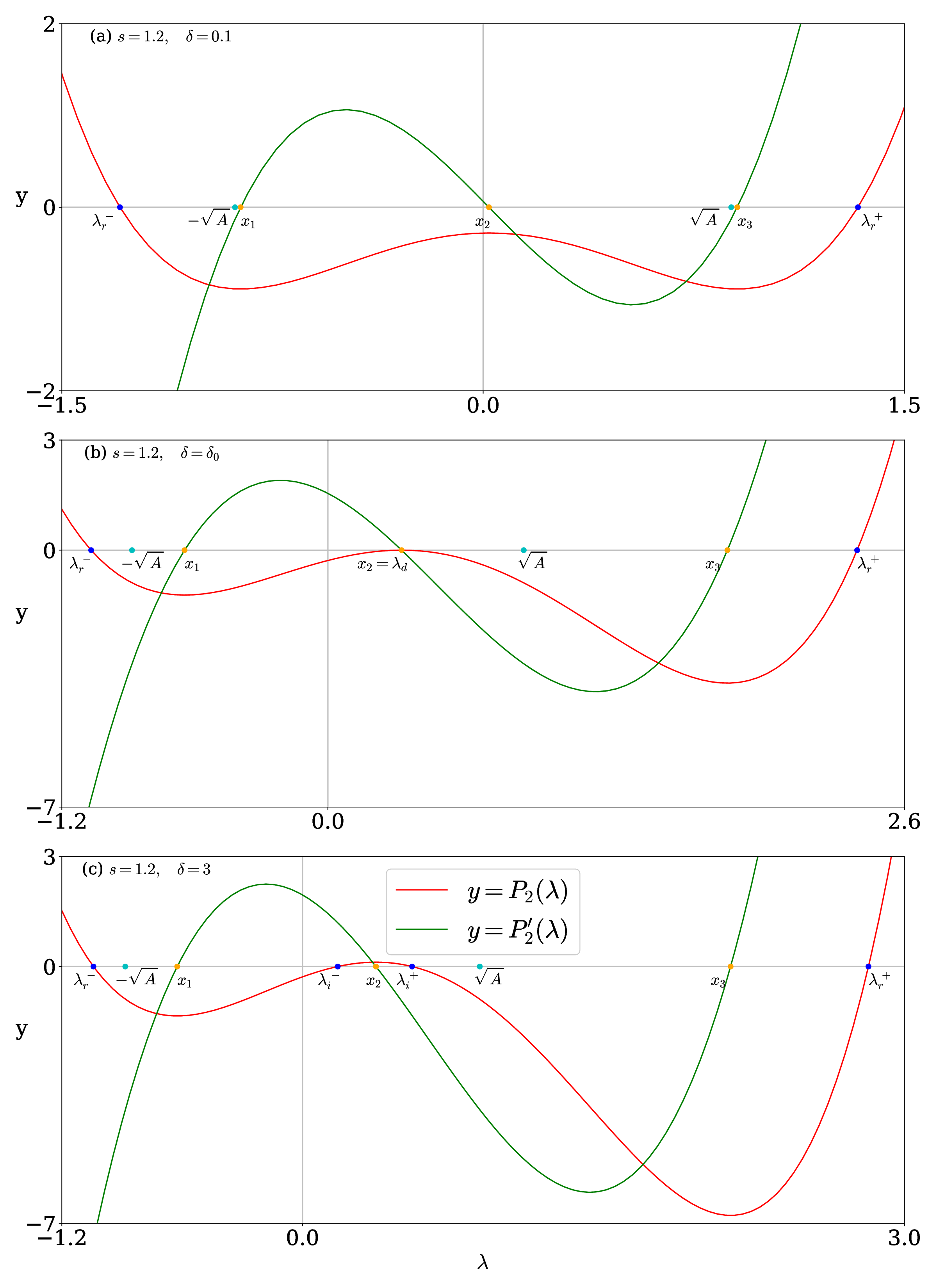}
\caption{The polynomials $P_2(\lambda)$ and $P_2'(\lambda)$ for $s=1.2$ and (a) $\delta<\delta_0$, (b) $\delta=\delta_0$ and (c) $\delta>\delta_0$. The real roots of $P_2(\lambda)$ are denoted by $\lambda_r^\pm$ and $\lambda_i^\pm$, and $\lambda_d$ is its double root at $\delta=\delta_0$. The roots of $P_2'(\lambda)$ are denoted by $x_1<0<x_2<x_3$. Note that the scales are different between graphs.}
\label{fig:polynomials}
\end{figure}

The double root of $P_2(\lambda)$ for $\delta=\delta_0$, which we will denote by $\lambda_d>0$, satisfies $\lambda_d<p_+<\lambda_r^+$, since $\lambda_r^+$ is the only root of $P_2(\lambda)$ in $(p_+,+\infty)$. Therefore, the double root is always less than $\lambda_r^+$. Figure \ref{fig:polynomials} shows the graph of the polynomials $P_2(\lambda)$ and $P'_2(\lambda)$ for $0<\delta<\delta_0$, $\delta=\delta_0$ and $\delta>\delta_0$ when $s=1.2$. 

Considering the real roots $\lambda=\lambda(\delta)=\lambda_r^\pm(\delta)$ of $P_2(\lambda)$ as functions of $\delta$, these are differentiable functions of $\delta$ since they are simple \cite[Theorem 6.3.12]{HJ2012}. Therefore, $\lambda(\delta)$ satisfies $P_2(\lambda(\delta))=0$ for all $\delta$. This implies that $\frac{d}{d\delta}P_2(\lambda(\delta))=0$ for all $\delta$. Note that because the coefficients of $P_2$ depend on $\delta$, then $\frac{d}{d\delta}P_2(\lambda(\delta))\not= P_2'(\lambda(\delta))\lambda'(\delta)$. Actually, differentiating the polynomial equation $P_2(\lambda(\delta))=0$ with respect to $\delta$, yields
\begin{equation}\label{eq:eigsderiv}
P_2'(\lambda(\delta))\frac{d}{d\delta}\lambda(\delta)=\frac{1}{s}\cdot\lambda(\delta)(\lambda(\delta)-\sqrt{A})(\lambda(\delta)+\sqrt{A})\ .
\end{equation}
For $\lambda=\lambda_r^-$, we have that $P_2'(\lambda_r^-)<0$ since $\lambda_r^-$ is less than the negative root of $P_2'(\lambda)$ (see also the Gauss--Lucas theorem). If $\lambda_r^-<-\sqrt{A}$, which is true for at least sufficiently small values of $\delta$, then $\frac{d}{d\delta}\lambda_r^-(\delta)>0$. The value $\pm \sqrt{A}$ cannot be a zero of $P_2(\lambda)$, since
$P_2(\pm\sqrt{A})=-\frac{1}{s(s-u_0)}\not=0$ for all $\delta>0$. Therefore, by the continuous dependence of the eigenvalues on $\delta$, we conclude $\lambda_r^-<-\sqrt{A}$ and $\frac{d}{d\delta}\lambda_r^-(\delta)>0$ for all $\delta>0$.

For $\lambda=\lambda_r^+$, we have that $P_2'(\lambda_r^+)>0$ since $\lambda_r^+$ is greater than the positive roots of $P_2'(\lambda)$, and $\lambda_r^+>\sqrt{A}$ for small values of $\delta$ and will remain as such since $\sqrt{A}$ cannot be a root of $P_2(\lambda)$. Thus, $\frac{d}{d\delta}\lambda_r^+(\delta)>0$, and $\lambda_r^+>\sqrt{A}$ for all $\delta>0$.

Note that $\lambda_r^+(0) + \lambda_r^-(0)=0$ and since $\frac{d}{d\delta}(\lambda_r^+(\delta) + \lambda_r^-(\delta))>0$, we have that 
$$\lambda_r^+(\delta) + \lambda_r^-(\delta)>0\ ,$$
for all $\delta>0$.

For $\delta>\delta_0$, the two real roots $\lambda_i^+$ and $\lambda_i^-$ of $P_2(\lambda)$ emerging from the double root $\lambda_d$  will be located left and right to the local maximum which is at $x_2$ (the second root of $P_2'(\lambda)$). Thus $P_2'(\lambda_i^-)>0$ and $P_2'(\lambda_i^+)<0$. Since $0<\lambda_i^-<\lambda_i^+<\sqrt{A}$, as it is shown below, the right-hand side of (\ref{eq:eigsderiv}) is negative for these roots, and we have that 
$$\frac{d}{d\delta}\lambda_i^-(\delta)<0\quad\text{and}\quad \frac{d}{d\delta}\lambda_i^+(\delta)>0\quad \text{for}\quad \delta>\delta_0\ ,$$
indicating that one moves to the left and the other to the right (see Figure \ref{fig:eigenvalues}).

Also, by writing the polynomial $P_2(\lambda)$ for $\delta>\delta_0$ in the form
$$P_2(\lambda)=(\lambda-\lambda_r^-)(\lambda-\lambda_r^+)(\lambda-\lambda_i^-)(\lambda-\lambda_i^+)\ ,$$
and equating its coefficients using the formula (\ref{eq:P2poly}) we obtain that
\begin{equation}\label{eq:maxeig}
\lambda_i^-=\frac{\delta}{s}-(\lambda_r^++\lambda_r^-+\lambda_i^+)<\frac{\delta}{s}\ ,
\end{equation}
since $\lambda_r^++\lambda_r^-+\lambda_i^+>0$. Similarly, $\lambda_i^+<\delta/s$.

The complex eigenvalues $\lambda=\lambda_i^\pm(\delta)$ emerging from the imaginary eigenvalues of the unperturbed system for $0<\delta<\delta_0$ are differentiable in $(0,\delta_0)$ and continuous for $\delta>0$. The differentiability breaks down when these collapse on the double positive eigenvalue $\lambda_d$. The double eigenvalue $\lambda_d=\lambda(\delta_0)$ lies in the interval $(0,\sqrt{A})$, as it is shown below.

Assume that for $0<\delta<\delta_0$ the complex eigenvalues are $\lambda_i^+(\delta)=z$ and $\lambda_i^-(\delta)=\bar{z}$. Then 
$$P_2(\lambda)=(\lambda-\lambda_r^+)(\lambda-\lambda_r^-)(\lambda-z)(\lambda-\bar{z})\ .$$
Executing the multiplications in the last expression, and equating the coefficients of $\lambda^3$ and of $\lambda^0$ terms with the corresponding terms of the formula (\ref{eq:P2poly}), we obtain
\begin{equation}\label{eq:reofz}
2 Re(z)+(\lambda_r^+ + \lambda_r^-)=\frac{\delta}{s}\ ,
\end{equation}
and 
\begin{equation}\label{eq:reofz2}
\lambda_r^+ \lambda_r^- |z|^2=-q\ ,
\end{equation}
respectively.

In order to treat the cases $\delta<\delta_0$ and $\delta\geq \delta_0$ simultaneously, we consider the symmetric functions
$$X(\delta)=\frac{\lambda_i^+(\delta)+\lambda_i^-(\delta)}{2}\quad\text{and}\quad m(\delta)=\lambda_i^+(\delta)\,\lambda_i^-(\delta)\ ,\qquad \delta\geq 0\ .$$
For $0<\delta<\delta_0$ we have $X=Re(z)$ and $m=|z|^2$, while $X(0)=0$ and $m(0)=\rho_i^2$. Since $\lambda_r^\pm$ are simple roots of $P_2(\lambda)$, the quadratic polynomial
$$\lambda^2-2X\lambda+m=\frac{P_2(\lambda)}{(\lambda-\lambda_r^+)(\lambda-\lambda_r^-)}$$
has coefficients that depend smoothly on $\delta>0$, also at $\delta=\delta_0$ where $\lambda_i^\pm$ are not differentiable. Denote $\sigma=\lambda_r^++\lambda_r^-$ and $p=\lambda_r^+\lambda_r^-<0$. Equating the coefficients of $P_2(\lambda)=(\lambda^2-\sigma\lambda+p)(\lambda^2-2X\lambda+m)$ with those of (\ref{eq:P2poly}) yields
\begin{equation}\label{eq:coefsys}
2X+\sigma=\frac{\delta}{s}\ ,\quad p+m+2X\sigma=-2A\ ,\quad 2Xp+\sigma m=-\frac{A\delta}{s}\ ,\quad pm=-q\ ,
\end{equation}
which generalize (\ref{eq:reofz}) and (\ref{eq:reofz2}). The last equation implies $m=-q/p>0$. Substituting $\sigma=\delta/s-2X$ and $p=-q/m$ into the third equation of (\ref{eq:coefsys}) we obtain
\begin{equation}\label{eq:Xofm}
X=\frac{\delta}{2s}\,\frac{m(m+A)}{m^2+q}\ ,
\end{equation}
and then the second equation of (\ref{eq:coefsys}) gives
\begin{equation}\label{eq:X2ofm}
4X^2(q-Am)=(q-2Am-m^2)(m+A)\ .
\end{equation}
Recall that $A^2+q=(1+\eta_0)/s^2$, and thus $\rho_i^2=-A+\sqrt{A^2+q}$ is the positive root of $q-2Am-m^2$. Moreover, $\rho_i^2<q/A$ since $A\rho_i^2<q$ is equivalent to $A\sqrt{A^2+q}<A^2+q$.

Since $m>0$, (\ref{eq:Xofm}) implies $X>0$ for $\delta>0$. Then by (\ref{eq:X2ofm}), for $m>0$ the quantities $q-Am$ and $q-2Am-m^2$ have the same sign, and thus $m(\delta)\in(0,\rho_i^2)\cup(q/A,\infty)$ for all $\delta>0$. Since $m(\delta)$ is continuous with $m(0)=\rho_i^2<q/A$, we conclude that $0<m(\delta)<\rho_i^2$ for all $\delta>0$. As $m<q/A$, we have $m(m+A)<m^2+q$, and (\ref{eq:Xofm}) gives
\begin{equation}\label{eq:mineig}
0<X(\delta)<\frac{\delta}{2s}\quad\text{for all}\quad \delta>0,\ s>1\ .
\end{equation}
In particular, $0<Re(\lambda_i^\pm)<\delta/2s$ for $0<\delta<\delta_0$ and $\lambda_d>0$. Moreover, for $\delta>\delta_0$ the roots $\lambda_i^\pm$ are real with positive sum $2X$ and positive product $m$, and thus $0<\lambda_i^-<\lambda_i^+$.
From (\ref{eq:inflpts}) we also observe that the inflection point $p_+>\delta/2s$. Thus, we have that 
$$0<Re(\lambda_i^\pm)<p_+<\lambda_r^+ \quad\text{for}\quad \delta>0 \ .$$

Combining (\ref{eq:Xofm}) and (\ref{eq:X2ofm}) we can express $\delta$ as a function of $m$:
\begin{equation}\label{eq:deltaofm}
\frac{\delta^2}{s^2}=R(m)\doteq\frac{(q-2Am-m^2)(m^2+q)^2}{(q-Am)\,m^2(m+A)}\ ,\qquad 0<m<\rho_i^2\ .
\end{equation}
A direct computation gives
$$R'(m)=-\frac{(m^2+q)\,K(m)}{m^3(m+A)^2(q-Am)^2}\ ,$$
where, denoting $r=q-2Am-m^2>0$,
$$K(m)=(2A+3m)\,r^3+7m(m+A)^2r^2+8m^2(m+A)^3r+4m^3(m+A)^4>0\ .$$
Thus, $R$ is strictly decreasing in $(0,\rho_i^2)$ with $R(\rho_i^2)=0$ and $\lim_{m\to0^+}R(m)=+\infty$. Consequently, $\delta\mapsto m(\delta)$ is a strictly decreasing bijection of $(0,\infty)$ onto $(0,\rho_i^2)$, and $m'(\delta)=2\delta/(s^2R'(m))<0$.

For $0<\delta<\delta_0$ we have $(Im(z))^2=m-X^2$, and from (\ref{eq:X2ofm})
\begin{equation}\label{eq:imofm}
m-X^2=\frac{C(m)}{4(q-Am)}\ ,\quad\text{where}\quad C(m)=m^3-Am^2+(2A^2+3q)m-Aq\ .
\end{equation}
The discriminant of $C'(m)=3m^2-2Am+2A^2+3q$ is $-20A^2-36q<0$, and thus $C$ is strictly increasing. Moreover, $C(0)=-Aq<0$ and $C(\rho_i^2)=4\rho_i^2(q-A\rho_i^2)>0$. Hence, $C$ has a unique root $m_d\in(0,\rho_i^2)$, and $\lambda_i^\pm$ are non-real if and only if $m(\delta)>m_d$, i.e. if and only if $0<\delta<s\sqrt{R(m_d)}$. Therefore, $\delta_0=s\sqrt{R(m_d)}$ and at $\delta=\delta_0$ we have $X^2=m_d$, so that $\lambda_d=\lambda(\delta_0)=\sqrt{m_d}$. Since $C(A)=2A(A^2+q)>0$ and $C$ is strictly increasing, we have $m_d<A$, and thus $\lambda_d<\sqrt{A}$. Since $\sqrt{A}$ is not a root of $P_2(\lambda)$ for any $\delta>0$, by continuity we also have $\lambda_i^+<\sqrt{A}$ for all $\delta>\delta_0$. Using (\ref{eq:Xofm}) we also obtain the alternative formula
$$\delta_0=\frac{2s(m_d^2+q)}{\sqrt{m_d}\,(m_d+A)}\ .$$
Furthermore, for $m_d<m<\rho_i^2$ we have
$$\frac{d}{dm}\left(m-X^2\right)=\frac{C'(m)(q-Am)+AC(m)}{4(q-Am)^2}>0\ ,$$
since $C'(m)>0$, $q-Am>0$ and $C(m)>0$ in this interval. Since $m'(\delta)<0$, we conclude that
$$\frac{d}{d\delta}(Im(z))^2=\frac{d}{dm}\left(m-X^2\right)m'(\delta)<0\quad\text{for}\quad 0<\delta<\delta_0\ .$$
Let $z=\lambda_i^+$ be the eigenvalue emerging from $+\I\rho_i$. Then $Im(\lambda_i^+)$ is continuous in $\delta$ with $Im(\lambda_i^+)=\rho_i>0$ at $\delta=0$, and it cannot vanish for $0<\delta<\delta_0$ since $\lambda_i^\pm$ are non-real there. Hence, $Im(\lambda_i^+)=\sqrt{m-X^2}>0$ and
$$\frac{d}{d\delta}Im(\lambda_i^+)=\frac{1}{2\,Im(\lambda_i^+)}\,\frac{d}{d\delta}(Im(z))^2<0\quad\text{for}\quad 0<\delta<\delta_0\ .$$
Since $\lambda_i^-=\overline{\lambda_i^+}$, we have $Im(\lambda_i^-)=-Im(\lambda_i^+)$, which is strictly increasing in $(0,\delta_0)$. Therefore, for all $s>1$, the complex eigenvalue emerging from $+\I\rho_i$ has strictly decreasing imaginary part, while the one emerging from $-\I\rho_i$ has strictly increasing imaginary part. Moreover, since $C(m_d)=0$, both imaginary parts tend to $0$ as $\delta\to\delta_0$, and the two eigenvalues collide on the double eigenvalue $\lambda_d$.

We note that the real part $Re(\lambda_i^\pm)$ is strictly increasing for sufficiently small $\delta>0$, in agreement with the asymptotic formula $\lambda_i^\pm\sim\frac{\delta}{4s}\pm\I\rho_i$, but it is not monotone in $(0,\delta_0)$ for all $s>1$. For example, for $s=2$ it attains its maximum at $\delta\approx 8.75<\delta_0\approx 14.54$.

Moreover, using an expansion of the form $\lambda=\lambda_0+\delta \lambda_1+\cdots$, we can find the asymptotic values of the eigenvalues as $\delta\to\infty$ to be:
$$ \begin{aligned}
\lambda_i^-&\sim\frac{1 + \eta_0 - (s- u_0)^2}{\delta (s - u_0)}+O(1/\delta^2)\ , \quad \lambda_i^+\sim \sqrt{\frac{s - u_0}{s}} -\frac{1 + \eta_0}{2 \delta (s - u_0)}+O(1/\delta^2)\ ,\\
\lambda_r^-&\sim -\sqrt{\frac{s - u_0}{s}} -\frac{1 + \eta_0 }{2 \delta (s - u_0)}+O(1/\delta^2)\ ,\quad 
\lambda_r^+\sim \frac{\delta}{s}+O(1)\ .
\end{aligned}
$$

Thus, we have that
$$\lim_{\delta\to\infty} \lambda_r^-=-\sqrt{A},\quad \lim_{\delta\to\infty} \lambda_r^+=\infty\ ,$$
and for $\delta>\delta_0$
\begin{multline*}-\rho_r<\lambda_r^-<-\sqrt{A}<0=\lim_{\delta\to \infty} \lambda_i^-<\lambda_i^-< \\ \lambda_i^+<\lim_{\delta\to +\infty} \lambda_i^+=\sqrt{A}<\rho_r<\lambda_r^+<\lim_{\delta\to\infty}\lambda_r^+=+\infty\ .
\end{multline*}
For $0<\delta<\delta_0$ the same ordering holds with $\lambda_i^\pm$ replaced by $Re(\lambda_i^\pm)$ only for moderate values of $s$. Numerically, we found that $Re(\lambda_i^\pm)$ exceeds $\sqrt{A}$ in a subinterval of $(0,\delta_0)$ when $s>4.9796$.

\bibliographystyle{plain}
\bibliography{biblio}

\end{document}